\documentclass[11pt, a4paper]{amsart}
\usepackage{fullpage} 
\usepackage{amsmath}
\usepackage{amsthm}
\usepackage{amssymb}
\usepackage[only,mapsfrom]{stmaryrd}
\usepackage{graphicx}
\usepackage{cancel} 
\usepackage{color} 
\usepackage[all]{xy} 
\usepackage{tikz-cd}
\usepackage{tikz}
\usetikzlibrary{arrows}
\usepackage{nicefrac}
\usepackage{enumitem}

\newtheorem{thm}{Theorem}[section]
\newtheorem{lem}[thm]{Lemma}
\newtheorem{prop}[thm]{Proposition}
\newtheorem{cor}[thm]{Corollary}

\newtheorem*{thm*}{Theorem}
\newtheorem*{conj*}{Conjecture}
\newtheorem*{cor*}{Corollary}
\newtheorem*{ques*}{Question}

\newtheorem*{namedthm}{\namedthmname}

\theoremstyle{definition}
\newtheorem{rem}[thm]{Remark}
\newtheorem{defn}[thm]{Definition} % definition numbers are dependent on theorem numbers
\newtheorem*{rem*}{Remark}

\DeclareMathOperator{\Lk}{Lk}
\DeclareMathOperator{\St}{St}

 \newcommand{\BR}{{\mathbb {R}}}

 \newcommand{\BZ}{{\mathbb {Z}}}
 \newcommand{\CB}{{\mathcal {B}}}
\newcommand{\CC}{{\mathcal {C}}}
 \newcommand{\CF}{{\mathcal {F}}}
\newcommand{\CG}{{\mathcal {G}}} \newcommand{\CH}{{\mathcal {H}}}
\newcommand{\CI}{{\mathcal {I}}} 
 
 \newcommand{\CN}{{\mathcal {N}}}

 \newcommand{\CT}{{\mathcal {T}}}

 \newcommand{\Aut}{{\mathrm{Aut}}}

\newcommand{\Hom}{{\mathrm{Hom}}}

\newcommand{\id}{{{\rm id}}}

\def\Proof{\noindent{\bf Proof.}\quad}
\def\sProof{\noindent{\bf Sketch of Proof. }\quad}
\def\qed{\hfill$\square$\smallskip}

\DeclareMathOperator{\Fix}{Fix}

\DeclareMathOperator{\Map}{Map}
\DeclareMathOperator{\RMap}{RMap}
\DeclareMathOperator{\PMap}{PMap}

\newenvironment{inlinecond}[1]{%
\begin{enumerate}%
\renewcommand*{\theenumi}{#1}%
\renewcommand*{\labelenumi}{#1}%
\item%
}{%
\end{enumerate}%
}

\title{Homological stability for the ribbon Higman--Thompson groups}

\author{Rachel Skipper}   
\address{Department of Mathematics, University of Utah, 155 S 1400 E,
Salt Lake City, UT 84112, USA} 
\email{rachel.skipper@utah.edu}

\author{Xiaolei Wu}
\address{Shanghai Center for Mathematical Sciences, Jiangwan Campus, Fudan University, No.2005 Songhu Road, Shanghai, 200438, P.R. China}
\email{xiaoleiwu@fudan.edu.cn}

\subjclass[2010]{20F36, 57M07, 19D23, 20J05}

\keywords{Braided Higman--Thompson groups, ribbon Higman--Thompson groups, asymptotic mapping class groups, big mapping class groups, finiteness property, homological stability. }

\begin{document}

\begin{abstract}
We generalize the notion of asymptotic mapping class groups and allow them to surject to the Higman--Thompson groups, answering a question of Aramayona and Vlamis in the case of the Higman--Thompson groups. When the underlying surface is a disk, these new  asymptotic mapping class groups can be identified with the ribbon and oriented ribbon Higman--Thompson groups. We use this model to prove that the ribbon Higman--Thompson groups satisfy homological stability, providing the first homological stability result for dense subgroups of big mapping class groups. Our result can also be treated as an extension of Szymik--Wahl's work on homological stability for the Higman--Thompson groups  to the surface setting.
\end{abstract}

\maketitle

\section*{Introduction}

The family of Thompson's groups and the many groups in the extended Thompson family have long been studied for their many interesting properties. Thompson's group $F$ is the first example of a type $F_\infty$, torsion-free group with infinite cohomological dimension \cite{BG84} while Thompson's groups $T$ and $V$ provided the first examples of finitely presented simple groups with infinitely many elements. More recently the braided and labeled braided Higman--Thompson groups have garnered attention in part due to their connections with big mapping class groups \cite{Br06,Deh06,AC20,SW21}.  In particular,  Thumann constructed the ribbon version of Thompson's group $V$ and proved that it is of type $F_\infty$ \cite{Th17}. The authors studied the ribbon Higman--Thompson groups $RV_{d,r}$ and their oriented version $RV_{d,r}^+$ in \cite{SW21}. In fact, we identified them with the so-called labeled braided Higman-Thompson groups and proved that they are all of type $F_\infty$.

The homology of Thompson's groups has also been well-studied.  Brown and Geoghegan computed the homology of $F$ in \cite{BG84}; Ghys and Sergiescu calculated the homology of $T$ in \cite{GS87}. More recently Szymik and Wahl showed that  $V$ is acyclic \cite{SW19}, answering a question due to Brown \cite{Bro89}. One of the key ingredients for their proof was showing that the Higman--Thompson groups $V_{d,1}\hookrightarrow V_{d,2}\hookrightarrow \cdots \hookrightarrow V_{d,r} \hookrightarrow \cdots$ satisfy homological stability for any fixed $d$. Recall that a family of groups 
$G_1\hookrightarrow G_2 \hookrightarrow \cdots \hookrightarrow G_n \hookrightarrow \cdots$
is said to satisfy homological stability if the induced maps $H_i(G_n)\rightarrow H_i(G_{n+1})$ are isomorphisms for sufficiently large $n$.  Classical examples of families of groups which satisfy homological stability include symmetric groups \cite{Nak61}, general linear groups \cite{vdK80} and mapping class groups of surfaces \cite{Har85}.

In the present paper, we extend  Szymik and Wahl's work to the class of ribbon Higman--Thompson groups. To accomplish this, we first build a geometric model for the ribbon Higman--Thompson groups using Funar--Kapoudjian's asymptotic mapping class groups \cite{FK04}. These groups are defined using a rigid structure on a surface minus a Cantor set and they sit naturally inside the ambient big mapping class groups. More recently, Aramayona and Funar \cite{AF17} generalized the definition to surfaces with nonzero genus.  In fact, Aramayona and Funar showed that the half-twist version of their asymptotic mapping class group (cf. Definition \ref{groupB}) is dense in the big mapping class group \cite[Theorem 1.3]{AF17}. Another surprising result of Funar and  Neretin says that the half-twist asymptotic mapping class group of a closed surface minus a standard ternary Cantor set is in fact isomorphic to its smooth mapping class group    \cite[Corollary 2]{FN18}.  In \cite[Question 5.37]{AV20}, the  following question was raised by Aramayona and Vlamis.

\begin{ques*}
Are there other geometrically defined subgroups of $\text{Map}(\Sigma_g)$ which surject to other interesting classes of subgroups of homeomorphism group of the Cantor set, such as the Higman--Thompson groups, Neretin groups, etc?
\end{ques*}

We proceed to construct two new classes of asymptotic mapping class groups, one of which answers their question in the case of Higman--Thompson groups while the other family surjects to the symmetric Higman--Thompson groups $V_{d,r}(\BZ/2\BZ)$.

\begin{thm*}[\ref{thm-ses-asm}, \ref{thm-dense-asy}]
Let $\Sigma$ be any compact surface and $\CC$ be a Cantor set which lies in the interior of a disk in $\Sigma$. Then the mapping class group $\Map(\Sigma\setminus \CC)$ contains the following two families of dense subgroups: the asymptotic mapping class groups $\CB V_{d,r}(\Sigma)$, which surject to the Higman--Thompson groups $V_{d,r}$, and the half-twist asymptotic mapping class groups $\CH V_{d,r}(\Sigma)$, which surject to the symmetric Higman--Thompson groups $V_{d,r}(\BZ/2\BZ)$.  
\end{thm*}

When $\Sigma$ is the disk, we identify $\CH V_{d,r}(\Sigma)$ with the ribbon Higman--Thompson group $RV_{d,r}$ and $\CB V_{d,r}(\Sigma)$  with the oriented ribbon Higman--Thompson group $RV^+_{d,r}$ (cf. Theorem \ref{thm-iden-asym-ribb}). Using this geometric model for the  ribbon Higman--Thompson groups, we are able to prove the following.

\begin{thm*}[\ref{thm-hmg-stab-RB+}, \ref{thm-hmg-stab-RB}]
Suppose $d\geq 2$. Then the inclusion maps induce isomorphisms
$$ \iota_{R,d,r}: H_i(RV_{d,r},M) \to H_i(RV_{d,r+1},M)$$
in homology in all dimensions $i \geq 0$, for all $r \geq  1$ and for all $H_1(RV_{d,\infty})$-modules $M$. The same also holds for the oriented  ribbon Higman--Thompson groups $RV_{d,r}^+$.
\end{thm*}
\begin{rem*}
\begin{enumerate}
\item Here we restrict our main result to the constant coefficient $\BZ$ case. Nevertheless,
the theorem also holds for some general coefficients by applying \cite[Theorem A]{RWW17}. 

\item  The same method here can also be used to prove that the groups $\CB V_{d,r}(\Sigma)$ and $\CH V_{d,r}(\Sigma)$ satisfy homological stability. Still it seems difficult to make it work directly for braided Higman-Thompson groups as we are lacking of a good geometric model for them. Ideally, we would realize the braided Higman-Thompson groups as some sort of mapping class groups of the disk minus Cantor set. In fact, since  the braided Higman-Thompson groups are subgroups of the oriented  ribbon Higman--Thompson groups, we  already have a geometric model in some sense. But it is less clear how one can tell when an element of the asymptotic mapping class group lies in the braided  Higman--Thompson groups.

\end{enumerate}
\end{rem*}

To the best of our knowledge, this is the first homological stability result for dense subgroups of big mapping class groups although density will not play a role in our proof. Our proof uses a recent convenient framework given by Randal-Williams and Wahl \cite{RWW17}. The core of the proof is similar to \cite{SW19}, but with new technical difficulties arising from infinite type surface topology. In particular, we take advantage of what we call the ``mutual link trick" which we abstract from \cite{BFM+16} and which we expect to be useful in a number of settings.
We hope our result here can be  further used to calculate the homology of ribbon Higman--Thompson groups and shed light on the question whether braided $V$ is acyclic. In fact, our Proposition \ref{prop-Tn-inf} has already been used in \cite{PalmerWu22} to prove that the mapping class groups of disk minus Cantor set is acyclic. It is also worth mentioning that, the homology of the infinite genus version ribbon Thompson group  has been calculated rationally in \cite[Theorem 1.2]{FK09} and integrally in \cite[Theorem 1.16]{ABF+21}. It in fact has the same homology as the stable homology of mapping class groups.

\subsection*{Outline of paper}
In Section \ref{sec:connectivitytools}, we describe the connectivity tools that will be necessary for the remainder of the paper. In Section \ref{sec:braidedgroups}, we introduce the definition of the Higman--Thompson, ribbon Higman--Thompson, and oriented ribbon Higman--Thompson groups using paired forest diagrams to define the elements. In Section \ref{sec-gm-bht-group}, we generalize the notion of asymptotic mapping class groups and allow them to surject to the Higman--Thompson groups. And finally, in Section \ref{sec:homstab}, we prove homological stability for the ribbon Higman--Thompson groups and their oriented version.

\subsection*{Notation and convention.} All surfaces in this paper are assumed to be connected and orientable unless otherwise stated. Given a simplicial complex $X$ and a cell $\sigma \in X$, we denote the link of $\sigma$ in $X$ by $\Lk_X(\sigma)$ (resp. the star of $\sigma$ by $\St_X (\sigma)$). When the situation is clear, we quite often omit $X$ and simply denote the link by $\Lk (\sigma)$ and the star by $\St (\sigma)$. Recall that $X$ is called $n$-connected if its homotopy groups are trivial up to dimension $n$.   We also use the convention that  $(-1)$-connected means non-empty and that every space is $(-2)$-connected. In particular, the empty set is $(-2)$-connected. Finally, we adopt the convention that elements in groups are multiplied from left to right.

\subsection*{Acknowledgements.}  
The first part of this project was done while the first author was a visitor in the Unit\'{e} de math\'{e}matiques pures et appliqu\'{e}es at the ENS de Lyon and during a visit to the University of Bonn. She thanks them for their hospitality. She was also supported by the GIF, grant I-198-304.1-2015, ``Geometric exponents of random walks and intermediate growth groups" and NSF DMS--2005297 ``Group Actions on Trees and Boundaries of Trees".

Part of this work was done when the second author was a member of the  Hausdorff Center of Mathematics. At the time, he was  supported by Wolfgang L\"uck's ERC Advanced Grant “KL2MG-interactions”
(no. 662400) and the DFG Grant under Germany's Excellence Strategy - GZ 2047/1, Projekt-ID 390685813.

Part of this work was also done when both authors were visiting IMPAN at Warsaw during the Simons Semester ``Geometric and Analytic Group Theory" which was partially supported by the grant 346300 for IMPAN from the Simons Foundation and the matching 2015-2019 Polish MNiSW fund. We would also like to thank Kai-Uwe Bux for inviting us for a research visit at Bielefeld in May 2019 and many stimulating discussions. Special thanks go to Jonas Flechsig for his comments on preliminary versions of the paper. Furthermore, we want to thank Javier Aramayona and Stefan Witzel for discussions, Andrea Bianchi for comments and Matthew C. B. Zaremsky for some helpful communications and comments. The authors are also grateful to the referee for many helpful comments which improved the exposition.

\section{Connectivity Tools}\label{sec:connectivitytools}
In this section, we review some of the connectivity tools that we need for calculating the connectivity of our spaces. A good reference is \cite[Section 2]{HV17} although not all the tools we use can be found there.

\subsection{Complete join}\label{subsec-com-join} The complete join is useful tool introduced by Hatcher and Wahl in \cite[Section 3]{HW10} for proving connectivity results. We review the basics here.

\begin{defn}
A surjective simplicial map $\pi: Y\to X$ is called a \emph{complete join} if it satisfies the following properties:
\begin{enumerate} [label=(\arabic*)]
    \item $\pi$ is injective on individual simplices. 
    \item For each $p$-simplex $\sigma = \langle v_0,\cdots,v_p\rangle$ of $X$, $\pi^{-1}(\sigma)$ is the join $\pi^{-1}(v_0)\ast \pi^{-1}(v_1)\ast\cdots\ast \pi^{-1}(v_p)$.
 \end{enumerate}
\end{defn}

%Recall a complex is \emph{$d$-spherical} if it is of dimension $d$ and is $(d-1)$-connected.

\begin{defn}
A simplicial complex  $X$ is called weakly Cohen-Macaulay of dimension $n$ if $X$ is $(n-1)$-connected and the link of each $p$-simplex of $X$ is $(n-p-2)$-connected. We sometimes shorten  weakly Cohen-Macaulay to $wCM$. 
\end{defn}

The main result regarding complete join that we will use is the following.

\begin{prop}\cite[Propostion 3.5]{HW10} \label{prop-join-conn}
If $Y$ is a complete join complex over a  $wCM$ complex $X$ of dimension $n$, then $Y$ is also  $wCM$ of dimension $n$.
\end{prop}

\begin{rem}\label{rem-cjoin}
If $\pi: Y\to X$ is a complete join, then $X$ is a retract of $Y$. In fact, we can define a simplicial map $s:X\to Y$ such that $\pi\circ s = \id_X$ by sending a vertex $v\in X$ to any vertex in $\pi^{-1}(v)$ and then extending it to simplices. The fact that $s$ can be extended to simplices is granted by the condition that $\pi$ is a complete join.  In particular we can also conclude that if $Y$ is $n$-connected, so is $X$.
\end{rem}

\subsection{Bad simplices argument}\label{sub-bad-sim}
Let $(X,Y)$ be a pair of simplicial complexes. We want to relate the $n$-connectedness of $Y$ to the $n$-connectedness of $X$ via a so called bad simplices argument, see \cite[Section 2.1]{HV17} for more information. One identifies a set of simplices in $X \setminus Y$ as bad simplices, satisfying the following two conditions:
\begin{enumerate}  [label=(\roman*)]
    \item \label{bad-sim-1} Any simplex with no bad faces is in $Y$, where by a ``face" of a simplex we mean a subcomplex spanned by any nonempty subset of its vertices, proper or not.
    \item \label{bad-sim-2} If two faces of a simplex are both bad, then their join is also bad.
\end{enumerate}

We call simplices with no bad faces good simplices. Bad simplices may have good faces or faces which are neither good nor bad. If $\sigma$ is a bad simplex, we say a simplex $\tau$ in $\Lk (\sigma)$ is good for $\sigma$ if any bad face of $\tau\ast \sigma$ is contained in $\sigma$. The simplices which are good for $\sigma$ form a subcomplex of $\Lk(\sigma)$ which we denote by $GL_{\sigma}$ and call the good link of $\sigma$.

\begin{prop} \cite[Proposition 2.1]{HV17}\label{prop-bad-sim}
Let $X, Y$ and $GL_\sigma$ be as above. Suppose that for some integer $n\geq 0$ the subcomplex $GL_\sigma$ of $X$ is $(n-\text{dim}(\sigma)-1)$-connected for all bad simplices $\sigma$. Then the pair $(X,Y)$ is $n$-connected, i.e. $\pi_i(X,Y) = 0$ for all $i\leq n$. 
\end{prop}

We can apply the proposition in the following way.

\begin{thm} \cite[Corollary 2.2]{HV17}\label{thm--bad-sim}
Let $Y$ be a subcomplex of a simplicial complex $X$ and suppose the space $X\setminus Y$ has a set of bad simplices satisfying \ref{bad-sim-1} and \ref{bad-sim-2} above, then:
\begin{enumerate}  [label=(\arabic*)]
    \item If $X$ is $n$-connected and $GL_\sigma$ is $(n-\text{dim}(\sigma))$-connected for all bad simplices $\sigma$, then $Y$ is $n$-connected.
    
    \item  If $Y$ is $n$-connected and $GL_\sigma$ is $(n-\text{dim}(\sigma)-1)$-connected for all bad simplices $\sigma$, then $X$ is $n$-connected. 
    
\end{enumerate}
\end{thm}

\subsection{The mutual link trick}
In the proof of \cite[Theorem 3.10]{BFM+16}, there is a beautiful argument for resolving intersections of arcs inspired by Hatcher's flow argument \cite{Hat91}. They attributed the idea to Andrew Putman. Recall Hatcher's flow argument allows one to ``flow" a complex to its subcomplex. But in the process, one can only ``flow" a vertex to a new one in its link. The mutual link trick will allow one to ``flow" a vertex to a new one not in its link provided ``the mutual link" is sufficiently connected.

To apply the  mutual link trick, we first need a lemma that allows us to homotope a simplicial map to a simplexwise injective one \cite[Lemma 3.9]{BFM+16}.  Recall a simplicial map is called \emph{simplexwise injective} if its
restriction to any simplex is injective. See also \cite[Section 2.1]{GRW18} for more information.

\begin{lem}
\label{lem:injectifying}
Let $Y$ be a compact $m$-dimensional combinatorial manifold.  Let $X$ be a
simplicial complex and assume that the link of every $p$-simplex in
$X$ is $(m-p-2)$-connected.  Let $\psi \colon Y \to X$ be a
simplicial map whose restriction to $\partial Y$ is simplexwise
injective.  Then after possibly subdividing the simplicial structure of $Y$, $\psi$ is
homotopic relative $\partial Y$ to a simplexwise injective map.
\end{lem}

Note that as discussed in \cite[Lemma 5.19]{GLU20}, there is a mistake in the connectivity bound given in \cite{BFM+16} that has been corrected later in \cite{BWZ21}. 

\begin{lem}[The mutual link trick]\label{lemma-replace-trick}
Let $Y$ be a closed $m$-dimensional combinatorial manifold and $f: Y\to X$ be a simplexwise injective simplicial map. Let $y \in Y$ be a vertex and $f(y) = x$ for some $x\in X$. Suppose $x'$ is another vertex of $X$ satisfying the following condition.
\begin{enumerate}
    \item $f(\Lk_{Y}(y)) \leq \Lk_{X} (x')$,
    \item the mutual link $\Lk_{X} (x) \cap \Lk_{X}(x')$ is $(m-1)$-connected,
\end{enumerate}
Then we can define a new simplexwise injective map $g:Y\to X $ by sending $y$ to $x'$ and all the other vertices $y'$ to $f(y')$ such that $g$ is homotopic to $f$. 
\end{lem}
\Proof 
The conditions that $f$ is simplexwise injective and $f(\Lk_{Y}(y)) \leq \Lk_{X} (x')$ guarantee that the definition of $g$ can be extended over $Y$ and $g$ is again simplexwise injective. 

We need to prove $g$ is homotopic to $f$. The homotopy will be the identity outside $\St_{Y}(y)$. Note that since $f$ is simplexwise injective, we have $f(\Lk_Y(y))\leq \Lk_X(x)$. Together with Condition (1), we have $f(\Lk_Y(y))\leq \Lk_X(x)\cap \Lk_{X} (x')$. Since
$\Lk_{Y}(y)$ is an~$(m-1)$-sphere and $\Lk_{X} (x) \cap \Lk_{X}(x')$ is $(m-1)$-connected,  there exists
an~$m$-disk $B$ with~$\partial B=\Lk_{Y}(y)$ and a simplicial
map~$\varphi\colon B \to \Lk_X(x)\cap\Lk_X(x')$ so that $\varphi$
restricted to $\partial B$ coincides with $\psi$ restricted
to~$\Lk_{Y}(y)$.  Since the image of~$B$ under~$\varphi$ is
contained in $\St_X(x)$ which is contractible, we can homotope~$g$,
replacing~$g|_{\St_{Y}(y)}$ with $\varphi$.  Since the image
of~$B$ under~$f$ is also contained in $\Lk_X(x')$, we can
similarly homotope $f$, replacing~$f|_{\St_{Y}(y)}$
with~$\varphi$.  These both yield the same map, so  $g$ is homotopic to $f$.  \qed

\section{Higman--Thompson groups and their braided versions}\label{sec:braidedgroups}
In this section, we first give an introduction to the Higman--Thompson groups and then define their ribbon version. 

\subsection{Higman--Thompson groups} The Higman--Thompson groups were first introduced by Higman as a generalization of the groups \cite{Hi74} given earlier in handwritten, unpublished notes of Richard Thompson.
First let us recall the definition of the Higman--Thompson groups. Although there are a number of equivalent definitions of these groups, we will use the notion of paired forest diagrams. First we define a \emph{finite rooted $d$-ary tree} to be a finite tree such that every vertex has degree $d+1$ except the \emph{leaves} which have degree 1, and the \emph{root}, which has degree $d$ (or degree $1$ if the root is also a leaf). Usually we draw such trees with the root at the top and the nodes descending from it down to the leaves. A vertex $v$ of the tree along with its $d$ adjacent descendants will be called a \emph{caret}. If the leaves of a caret in the tree are leaves of the tree, we will call the caret \emph{elementary}.  A collection of $r$ many $d$-ary trees will be called a $\emph{$(d,r)$-forest}$. When $d$ is clear from the context, we may just call it an $r$-forest.

Define a \emph{paired $(d,r)$-forest diagram} to be a triple $(F_-,\rho,F_+)$
consisting of two $(d,r)$-forests $F_-$ and $F_+$ both with $l$ leaves for some $l$, and a permutation $\rho \in S_l$, the symmetric group on $l$ elements.  We label the leaves
of~$F_-$ with $1,\dots,l$ from left to right, and for each $i$,
the $\rho(i)^{\text{th}}$ leaf of~$F_+$ is labeled $i$.

\begin{figure}[h]\label{fig:ele-V}
\centering
\begin{tikzpicture}[line width=1pt, scale=0.5]
\begin{scope}[xshift= -4cm]
  \draw
   (-12,-2) -- (-9,0) -- (-6,-2)
   (-9,0) -- (-9,-2)
   (-11,-1.33) -- (-11,-2)
   (-11,-1.33) -- (-10,-2)
   (-4,-2) -- (-1,0) -- (2,-2)
   (-1,0)  -- (-1,-2);

\draw[dotted] (3.5,0)  -- (3.5,-3);

  \filldraw
  (-12,-2) circle (1.5pt)
  (-9,0) circle (1.5pt)
   (-6,-2) circle (1.5pt)
   (-9,-2) circle (1.5pt)
  (-11,-1.33)  circle (1.5pt)
   (-11,-2) circle (1.5pt)
   (-10,-2) circle (1.5pt)
    (-4,-2) circle (1.5pt)
   (-1,0) circle (1.5pt)
   (-1,-2) circle (1.5pt)
  % (-3.5,0) circle (1.5pt)
   (2,-2) circle (1.5pt);
  \node at  (-12,-2.5) {$1$};
  \node at (-11,-2.5) {$2$};
  \node at (-10,-2.5) {$3$};
  \node at (-9,-2.5) {$4$};
  \node at (-6,-2.5) {$5$};
  \node at (-4,-2.5) {$6$};
  \node at (-1,-2.5) {$7$};
  \node at (2,-2.5) {$8$};
\end{scope}

  \begin{scope}[ xshift=13cm]
   \draw
   (-12,-2) -- (-9,0) -- (-6,-2)
   (-9,0) -- (-9,-2)
   
  (-8,-2) -- (-9,-1.33) -- (-10,-2)
   
   (-4,-2) -- (-1,0) -- (2,-2)
   (-1,0)  -- (-1,-2);

  \filldraw
  (-12,-2) circle (1.5pt)
  (-9,0) circle (1.5pt)
   (-6,-2) circle (1.5pt)
   (-9,-2) circle (1.5pt)
  (-9,-1.33)  circle (1.5pt)
   (-10,-2) circle (1.5pt)
   (-8,-2) circle (1.5pt)
    (-4,-2) circle (1.5pt)
   (-1,0) circle (1.5pt)
   (2,-2) circle (1.5pt)
   (-1,-2) circle (1.5pt);
  \node at  (-12,-2.5) {$4$};
  \node at (-6,-2.5) {$6$};
  \node at (-10,-2.5) {$1$};
  \node at (-8,-2.5) {$3$};
  \node at (-9,-2.5) {$2$};
  \node at (-4,-2.5) {$7$};
  \node at (2,-2.5) {$5$};
  \node at (-1,-2.5) {$8$};
 \end{scope}

  \begin{scope}[ yshift=-4.5cm,xshift= -4cm]
    \draw
   (-12,-2) -- (-9,0) -- (-6,-2)
   (-9,0) -- (-9,-2)
  % (-11,-1.33) -- (-11,-2)
  % (-11,-1.33) -- (-10,-2)
   (-4,-2) -- (-1,0) -- (2,-2)
   (-1,0)  -- (-1,-2);

\draw[dotted] (3.5,0)  -- (3.5,-3);

  \filldraw
  (-12,-2) circle (1.5pt)
  (-9,0) circle (1.5pt)
   (-6,-2) circle (1.5pt)
   (-9,-2) circle (1.5pt)
  %(-11,-1.33)  circle (1.5pt)
  % (-11,-2) circle (1.5pt)
  % (-10,-2) circle (1.5pt)
    (-4,-2) circle (1.5pt)
   (-1,0) circle (1.5pt)
   (-1,-2) circle (1.5pt)
  % (-3.5,0) circle (1.5pt)
   (2,-2) circle (1.5pt);
  \node at  (-12,-2.5) {$1$};
 % \node at (-11,-2.5) {$2$};
 % \node at (-10,-2.5) {$3$};
  \node at (-9,-2.5) {$2$};
  \node at (-6,-2.5) {$3$};
  \node at (-4,-2.5) {$4$};
  \node at (-1,-2.5) {$5$};
  \node at (2,-2.5) {$6$};
  \end{scope}

  \begin{scope}[xshift=13cm, yshift=-4.5cm]
   \draw
   (-12,-2) -- (-9,0) -- (-6,-2)
   (-9,0) -- (-9,-2)
   
  %(-8,-2) -- (-9,-1.33) -- (-10,-2)
   
   (-4,-2) -- (-1,0) -- (2,-2)
   (-1,0)  -- (-1,-2);

  \filldraw
  (-12,-2) circle (1.5pt)
  (-9,0) circle (1.5pt)
   (-6,-2) circle (1.5pt)
   (-9,-2) circle (1.5pt)
 % (-9,-1.33)  circle (1.5pt)
  % (-10,-2) circle (1.5pt)
  % (-8,-2) circle (1.5pt)
    (-4,-2) circle (1.5pt)
   (-1,0) circle (1.5pt)
   (2,-2) circle (1.5pt)
   (-1,-2) circle (1.5pt);
  \node at  (-12,-2.5) {$2$};
  \node at (-6,-2.5) {$4$};
 % \node at (-10,-2.5) {$1$};
%  \node at (-8,-2.5) {$3$};
  \node at (-9,-2.5) {$1$};
  \node at (-4,-2.5) {$5$};
  \node at (2,-2.5) {$3$};
  \node at (-1,-2.5) {$6$};
  \end{scope}
\end{tikzpicture}

\caption{Reduction, of the top paired $(3,2)$-forest diagram to the bottom one.}
\label{fig:reduction_V}
\end{figure}
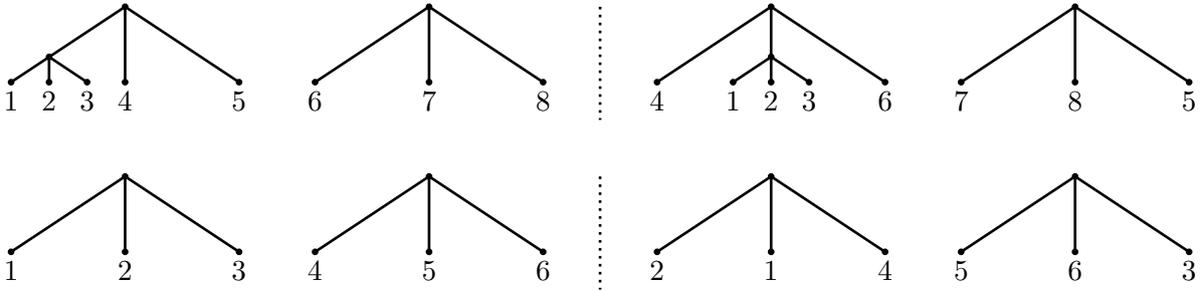

Define a \emph{reduction} of a paired $(d,r)$-forest diagram to be the following: Suppose there is an
elementary caret in~$F_-$ with leaves labeled by $i,\cdots,i+d-1$ from left to right, and an elementary caret in~$F_+$ with  leaves labeled by ~$i, \cdots,i+d-1$ from left to right.  Then we can ``reduce'' the diagram
by removing those carets, renumbering the leaves and replacing
$\rho$ with the permutation~$\rho'\in S_{l-d+1}$ that sends the new leaf
of~$F_-$ to the new leaf of~$F_+$, and otherwise behaves like~$\rho$.
The resulting paired forest diagram~$(F'_-,\rho',F'_+)$ is then said to
be obtained by \emph{reducing}~$(F_-,\rho, F_+)$. See Figure~\ref{fig:reduction_V} below for an idea of reduction
of paired $(3,2)$-forest diagrams. The reverse
operation to reduction is called \emph{expansion}, so $(F_-,\rho,F_+)$
is an expansion of $(F'_-,\rho',F'_+)$.  A paired forest diagram is
called \emph{reduced} if there is no reduction possible.  Define an equivalence relation on the set of paired $(d,r)$-forest diagrams by declaring two paired  forest diagrams to be equivalent if one can be reached by the other through a finite series of reductions and expansions.
Thus an equivalence class of paired forest diagrams consists of all diagrams
having a common reduced representative.  Such reduced representatives
are unique. 

There is a binary operation $\ast$ on the set of equivalence classes
of paired $(d,r)$-forest diagrams.  Let $\alpha =(F_-,\rho,F_+)$ and $\beta=(E_-,\xi,E_+)$
be reduced paired forest diagrams.  By applying repeated expansions to
$\alpha$ and $\beta$ we can find representatives $(F'_-,\rho',F'_+)$ and
$(E'_-,\xi',E'_+)$ of the equivalence classes of $\alpha$ and
$\beta$, respectively, such that~$F'_+ = E'_-$.  Then we
declare $ \alpha \ast  \beta $ to be $(F'_-,\rho'\xi',E'_+)$.  This operation is well defined on the equivalence classes and is a group operation.

\begin{defn}
    \label{def:Higman-V-F-T}
   The \emph{Higman--Thompson group} $V_{d,r}$ is the group of equivalence classes of paired
   $(d,r)$-forest diagrams with the multiplication~$\ast$.  
\end{defn}

The usual Thompson group $V$ is a special case of Higman--Thompson groups. In fact, $V= V_{2,1}$. 

\subsection{Ribbon Higman--Thompson groups}\label{subsect:braidedHTgroups} 

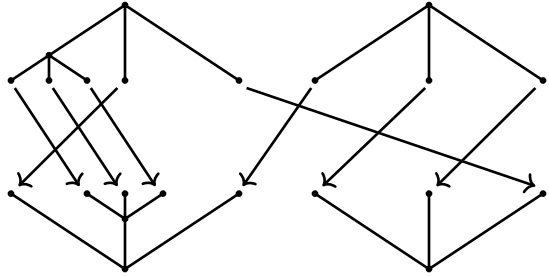
\begin{figure}[h]
\centering
\begin{tikzpicture} [line width=1pt, scale = 0.5]
 \begin{scope}[xshift =10cm]
 
 \draw
   (-12,-2) -- (-9,0) -- (-6,-2)
   (-9,0) -- (-9,-2)
   (-11,-1.33) -- (-11,-2)
   (-11,-1.33) -- (-10,-2)
   (-4,-2) -- (-1,0) -- (2,-2)
   (-1,0)  -- (-1,-2);
     %    \draw[dotted] (-1,0)  -- (-1,-3);
   \draw[->](-11.9,-2.2) -> (-10.2,-4.8);
      \draw[->](-10.9,-2.2) -> (-9.2,-4.8);
      \draw[->](-9.9,-2.2) -> (-8.2,-4.8);
            \draw[->](-9.2,-2.2) -> (-11.8,-4.8);
      
        \draw[->](-5.8,-2.2) -> (1.8,-4.8);
         \draw[->](-4.1,-2.2) -> (-5.9,-4.8);
        \draw[->](-1.1,-2.2) -> (-3.8,-4.8);
          \draw[->](1.8,-2.2) -> (-0.8,-4.8);
      
  \filldraw
  (-12,-2) circle (1.5pt)
  (-9,0) circle (1.5pt)
   (-6,-2) circle (1.5pt)
   (-9,-2) circle (1.5pt)
  (-11,-1.33)  circle (1.5pt)
   (-11,-2) circle (1.5pt)
   (-10,-2) circle (1.5pt)
    (-4,-2) circle (1.5pt)
   (-1,0) circle (1.5pt)
   (-1,-2) circle (1.5pt)
  % (-3.5,0) circle (1.5pt)
   (2,-2) circle (1.5pt);

\end{scope}

 \begin{scope}[ xshift = 10cm, yscale=-1,  yshift=7cm]
  \draw
   (-12,-2) -- (-9,0) -- (-6,-2)
   (-9,0) -- (-9,-2)
   
  (-8,-2) -- (-9,-1.33) -- (-10,-2)
   
   (-4,-2) -- (-1,0) -- (2,-2)
   (-1,0)  -- (-1,-2);

  \filldraw
  (-12,-2) circle (1.5pt)
  (-9,0) circle (1.5pt)
   (-6,-2) circle (1.5pt)
   (-9,-2) circle (1.5pt)
  (-9,-1.33)  circle (1.5pt)
   (-10,-2) circle (1.5pt)
   (-8,-2) circle (1.5pt)
    (-4,-2) circle (1.5pt)
   (-1,0) circle (1.5pt)
   (-1,-2) circle (1.5pt)
  % (-3.5,0) circle (1.5pt)
   (2,-2) circle (1.5pt);
 \end{scope}   
   
\end{tikzpicture}
\caption{An element of $V_{3,2}$.}
\label{fig:element_of_V}
\end{figure}

For convenience, we will think of 
the forest $F_+$ drawn beneath $F_-$ and upside down, i.e., with the
root at the bottom and the leaves at the top.  The permutation $\rho$
is then indicated by arrows pointing from the leaves of $F_-$ to the
corresponding paired leaves of $F_+$.  See Figure~\ref{fig:element_of_V}
for this visualization of (the unreduced representation of) the
element of $V_{3,2}$ from Figure~\ref{fig:reduction_V}.

Now in the ribbon version of the Higman--Thompson groups, the permutations of leaves are simply
replaced by ribbon braids which can twist between the leaves.

\begin{defn}\label{defn-rb}
Let  $\CI = \amalg_{i=1}^d I_i: [0,1] \times \{1,\cdots,l\} \to \BR^2 $ be an embedding which we refer to as the \emph{marked bands}. A \emph{ribbon braid} is a map  $R: ([0,1] \times \{0,1,\cdots,l\})  \times [0,1] \to  \BR^2$  such that for any $0\leq t\leq 1$, $R_t: [0,1] \times \{1,\cdots,l\} \to  \BR^2$ is an embedding,  $R_0 = \CI$, and there exists $\sigma\in S_l$ such that $R_1(t) \mid_{I_i} = I_{\sigma(i)}(t)$ or $R_1(t)|_{I_i} = I_{\sigma(i)}(1-t)$. The usual product of paths defines a group structure on the set of ribbon braids up to homotopy among ribbon braids. This group, denoted by $RB_l$, does not depend on the
choice of the marked bands and it is called the ribbon braid group  with $l$ bands. A ribbon braid is \emph{pure} if $\sigma$ is trivial and we define $PRB_l$ to be the \emph{pure ribbon braid group} with $l$ bands. If we further assume $R_1(t) \mid_{I_i} = I_{\sigma(i)}(t)$, this subgroup is called the \emph{oriented ribbon braid group} $RB_l^+$. Similarly, we have the \emph{oriented pure ribbon braid group} $PRB_l^+$. 
\end{defn}

\begin{rem}
Note that $RB_l \cong \BZ^l\rtimes B_l$ where the action of the braid group  $B_l$ with $l$ strings is induced by the symmetric group action on the coordinates of $\BZ^l$. In particular, for the pure ribbon braid group $PRB_l$, we have $PRB_l  \cong \BZ^l\times PB_l$, where $PB_l$ is the pure braid group with $l$ strings. Under this isomorphism, $RB_l^+ \cong (2\BZ)^l\rtimes B_l$ and $PRB^+_l\cong (2\BZ)^l\times PB_l$.
\end{rem}

\begin{defn}
    A \emph{ribbon braided paired $(d,r)$-forest diagram} is a triple $(F_-,\mathfrak{r},F_+)$
    consisting of two $(d,r)$-forests $F_-$ and $F_+$ both with $l$ leaves for some $l$
     and a ribbon braid~$\mathfrak{r} \in RB_l$ connecting the leaves of $F_-$ to the leaves of $F_+$.
\end{defn}

\begin{figure}[h]
\centering
\begin{tikzpicture}[line width=1pt,scale=1.2]
 \draw
   (0,0) -- (1,0)   (0,-2) -- (1,-2) 
  (0,0) to [out=-90, in=90] (1,-2)   (1,0) to [out=-90, in=90] (0,-2);
  
    \draw[white, line width=4pt]
    (0,-0.1) to [out=-90, in=90] (1,-1.9);
    
   \draw   (0,0) to [out=-90, in=90] (1,-2);  
    
  \draw[->]
   (3,-1) -> (4,-1);

 \begin{scope}[xshift=6.5cm]

   \draw (0,0) -- (0.3,0) (0.7,0) -- (1,0)   (0,-2) -- (0.3,-2) (0.7,-2) -- (1,-2);

    %  \draw   (1,0) to [out=-90, in=90] (0.6,-.7);
      
      \draw (0.3,0) to [out=-90, in=90] (0.1,-0.7)
    (1,0) to [out=-90, in=60] (0.6,-.7);
  
    \draw[white, line width=4pt]
 (0,-.1) to [out=-90, in=110] (0.4,-0.7)
 (.7,-.1) to [out=-90, in=80] (0.9,-.7) ;
 
\draw (0,0) to [out=-90, in=120] (0.4,-0.7)  (0.7,0) to [out=-90, in=90]  (0.9,-.7)  ;

 \draw   
 (0.9,-0.7) to [out=-120, in=90] (0.3,-2)
 (0.6,-0.7) to [out=-120, in=90] (0,-2);
    
   \draw[white, line width=4pt]
(0.1,-0.7) to [out=-90, in=90] (0.7,-2)
 (0.4,-0.7) to [out=-60, in=90] (1,-2) ;

     \draw   
 (0.1,-0.7) to [out=-90, in=90] (0.7,-2)
 (0.4,-0.7) to [out=-60, in=90] (1,-2);

 \end{scope}
\end{tikzpicture}
\caption{Splitting a ribbon into $2$ ribbons.}
\label{fig:Splitting_ribb}
\end{figure}
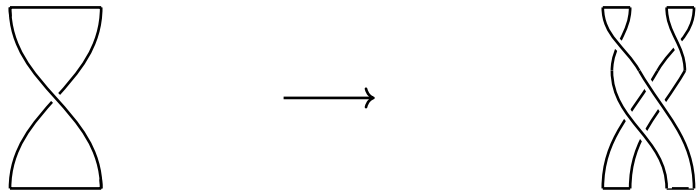

The expansion and reduction rules for the ribbon braids just come from the natural way of splitting a ribbon band into $d$ components and the inverse operation to this. See Figure \ref{fig:Splitting_ribb} for how to split a half twisted band when $d=2$. Note that not only are the two bands themselves twisted but the bands are also braided.  Everything else will be the same as in the braided case, so we omit the details here.  As usual, we define two ribbon braided paired forest diagrams to be equivalent if one is obtained
from the other by a sequence of reductions or expansions.  The
multiplication operation~$\ast$ on the equivalence classes is defined
the same way as for~$bV_{d,r}$. We direct the reader to \cite[Section 2]{SW21}

\begin{defn}
    The \emph{ribbon Higman--Thompson group} $RV_{d,r}$ (resp. the \emph{oriented ribbon Higman--Thompson group} $RV^+_{d,r}$) is the group of equivalence
    classes of (resp. oriented) ribbon braided paired $(d,r)$-forests diagrams with the
    multiplication~$\ast$.
\end{defn}

\section{Asymptotic mapping class groups related to the ribbon Higman--Thompson groups}\label{sec-gm-bht-group}
The purpose of this section is to  generalize the notion of asymptotic mapping class groups and allow them to surject to the Higman--Thompson groups.   In particular, we  will build a geometric model for the ribbon Higman--Thompson groups which will be crucial for proving homological stability in Section~\ref{sec:homstab}. Our construction is largely based on the ideas in  \cite[Section 2]{FK04} and \cite[Section 3]{AF17}.

\subsection{$d$-rigid structure} In this subsection, we generalize the notion of a rigid structure to that of a $d$-rigid structure. 

\begin{defn}
A {\em  $d$-leg pants} is a surface which is homeomorphic to a $(d+1)$-holed sphere.  
\end{defn}

\begin{figure}
\centering
\includegraphics[width=0.7\textwidth]{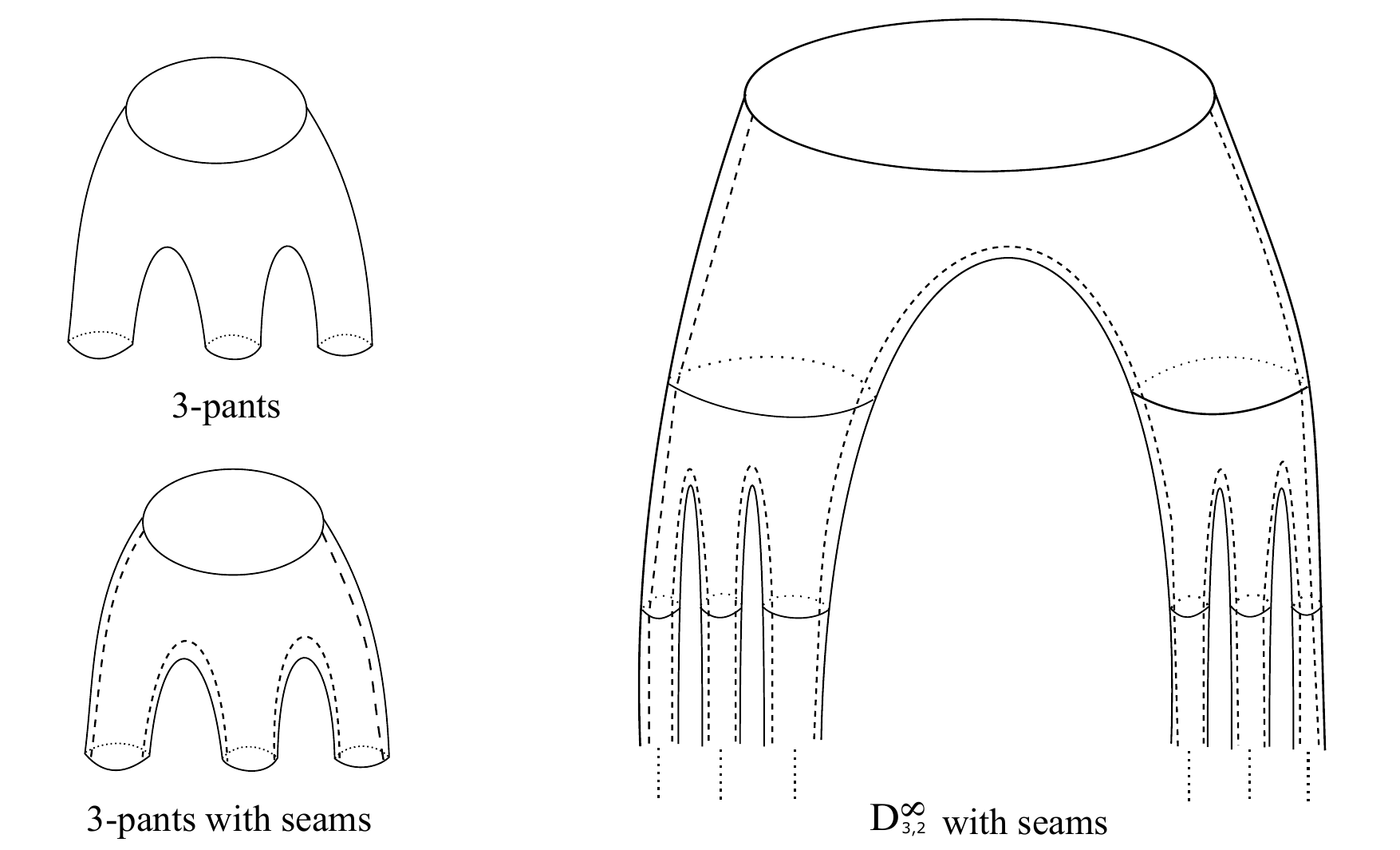}
\caption{$3$-leg pants and the surface $D_{3,2}^\infty$ with canonical seams}
\label{fig-inf-surf-seam}
\end{figure}

Recall that the usual pair of pants is a $2$-leg pants. We will draw a $d$-leg pants with one boundary component at the top. In this way, we can conveniently put a counter-clockwise total order on the boundary components, making the top component the minimal one. See Figure \ref{fig-inf-surf-seam} for an example of a $3$-leg pants.

We proceed to build some infinite type surfaces using some basic building blocks.

\begin{defn}\label{defn-inf-surf}
Let $\Sigma$ be an compact oriented surface. Call the boundary components of $\Sigma$ the \emph{based boundary components}. Then $\Sigma^{\infty}_{d,r}$ is the infinite surface, built up as an inductive limit of infinite surfaces $\Sigma_{d,r,m}$ with $m\geq 0$:
\begin{enumerate}
    \item $\Sigma_{d,r,0}$ is obtained from $\Sigma$ by deleting the interior of a  disk in $\Sigma$.  When $\Sigma$ is a disk $D$, we declare $D_{d,r,0} =\partial D$.
    \item  $\Sigma_{d,r,1}$ is obtained from $\Sigma_{d,r,0}$ attaching a copy of $r$-leg pants along the newly created boundary of $\Sigma_{d,r,0}$.
    \item For $m\geq 1$, $\Sigma_{d,r,m+1}$ is obtained from $\Sigma_{d,r,m}$ by gluing  a pair of $d$-leg pants to  every nonbased boundary circle of $\Sigma_{d,r,m}$ along the top boundary of the pants.
\end{enumerate}
The surface $\Sigma_{d,r,1}$ is called the \emph{base} of $\Sigma^{\infty}_{d,r}$ and the boundary components of  $\Sigma^{\infty}_{d,r}$ coming from the base are the \emph{based boundary components}. For each $m \geq 1$, the  nonbased boundary components of $\Sigma_{d,r,m}$ naturally embed in $\Sigma^{\infty}_{d,r}$ and we call these the \emph{admissible loops}. We call the admissible loops coming from $\Sigma_{d,r,1}$ the \emph{rooted loops}. The surface $\Sigma^{\infty}_{d,r}$ has a natural induced orientation.
\end{defn}

\begin{rem}
The two indices $d,r$ in the definition of $\Sigma^{\infty}_{d,r}$ will be used later to define the Higman--Thompson version of the asymptotic mapping class group (see Definition \ref{groupB}), where $d$ is related to the valance of the rooted trees and $r$ is the number of roots in the definition of the Higman--Thompson groups.  
\end{rem}

\begin{rem}
 To define our $d$-rigid structure, we do not really need $\Sigma_{d,r,0}$. But it will be  convenient  to have $\Sigma_{d,r,0}$   later  in Definition \ref{defn-forest} for defining the map from  $\Sigma^{\infty}_{d,r}$ to the tree $\CT_{d,r}$.
\end{rem}

\begin{defn}\label{defn-add-subsur}
 A compact subsurface $A \subset \Sigma^{\infty}_{d,r}$ is {\em admissible} if $\Sigma_{d,r,1} \subseteq A$ and all of its nonbased boundaries are admissible. The subsurfaces $\Sigma_{d,r,m}$ are called the \emph{standard admissible subsurfaces} of $\Sigma^{\infty}_{d,r}$. 
\end{defn}

\begin{figure}
\centering
\includegraphics[width=0.5\textwidth]{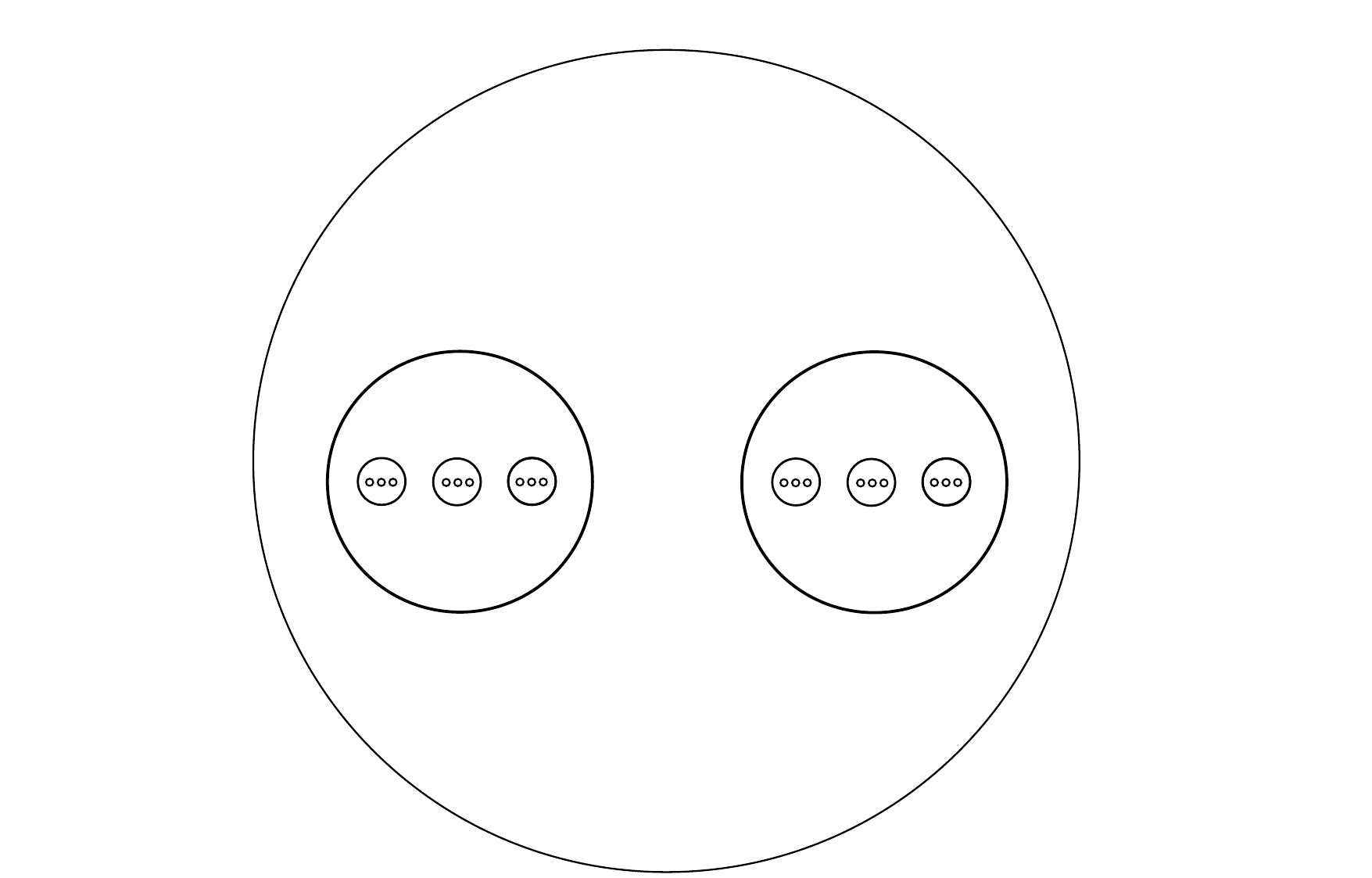}
\caption{Disk model for the surface $D_{3,2}^{\infty}$}
\label{disk-model-inf-surface}
\end{figure}

\begin{rem}\label{rem-emb-model-disk}
In the special case where the starting surface is a disk, we will use the notation $\Sigma=D$, $\Sigma_{d,r,m}=D_{d,r,m}$, and $\Sigma^{\infty}_{d,r}=D^{\infty}_{d,r}$. See Figure \ref{fig-inf-surf-seam} for a picture of the surface $D^{\infty}_{3,2}$. 
In this case, we can think of $D^{\infty}_{d,r}$ as a subsurface of a disk $D$. More specifically, let $D= \{(x,y)\mid x^2 + y^2 \leq 1\}$ and $x_i= \frac{2i-r-1}{r+1}$, $1\leq i\leq r$. We place $r$ disks with center at each $(x_i,0)$ of radius $r_0=\frac{1}{4(r+1)}$. Denote these disks by $D_1,\cdots D_r$. The complement of the interior of these $r$ disks in $D$ is homeomorphic to the $r$-leg pants $D_{d,r,1}$. Now for each disk  $D_i$, $1\leq i \leq r$, we can equally distribute $d$ points in the $x$-axis inside $D_i$ and place a disk with radius $\frac{r_0}{d^2}$ centered at each. We have the complement of the interior of these $d$ disks in $D_i$ are all  $d$-leg pants. We can continue the process inductively. At the end, the disks converge to a Cantor set which we denote by $\CC$. In particular $D^{\infty}_{d,r}$ is homeomorphic to $D \setminus \CC$.  We will refer to this as the {\em puncture model} for $D^{\infty}_{d,r}$. See Figure \ref{disk-model-inf-surface} for a picture of $D_{3,2}^\infty$ with this model.  The advantage of this model is we can view $D^{\infty}_{d,r}$ and all its admissible subsurfaces directly as a subsurfaces of $D$. 
\end{rem}

\begin{rem}\label{rem-emb-model-gen}
 Now $\Sigma^{\infty}_{d,r}$ can be obtained from $\Sigma$ by attaching a copy of $D^{\infty}_{d,r}$ to the nonbased boundary component of $\Sigma_{d,r,0}$. In particular, $\Sigma^{\infty}_{d,r}$ is obtained from  $\Sigma$ by deleting a copy of the Cantor set, and  any admissible subsurface of  $\Sigma^{\infty}_{d,r}$ can be viewed directly as a subsurface of $\Sigma$ using the puncture model. Recall that any two Cantor sets are homeomorphic, hence,  by the classification of infinite surfaces \cite[Theorem 2.2]{AV20}, we have $\Sigma_{d,r}^\infty$ is homeomorphic to $\Sigma\setminus \CC $ where $\CC$ is the standard ternary Cantor set sitting inside some disk in $\Sigma$ regardless of the choice of $d$ and $r$.
\end{rem}

\begin{defn}\noindent
\begin{enumerate}
\item A {\em suited $d$-pants decomposition} of the infinite surface $\Sigma^{\infty}_{d,r}$ is a maximal collection
of distinct nontrivial simple closed curves  in the interior of  $\Sigma^{\infty}_{d,r} \setminus \Sigma_{d,r,1}$ which are not isotopic to the boundary, pairwise disjoint and pairwise non-isotopic, with the additional property that the complementary regions in  $\Sigma^{\infty}_{d,r} \setminus \Sigma_{d,r,1}$ are  all $d$-leg pants.
\item A {\em $d$-rigid structure}  on $\Sigma^{\infty}_{d,r}$ consists of two
pieces of data:%\vspace{0.1cm}\\
\begin{itemize}
\item a suited $d$-pants
  decomposition, and
\item  a {\em $d$-prerigid structure}, i.e. a countable collection of disjoint
line segments embedded into $\Sigma^{\infty}_{d,r} \setminus \Sigma_{d,r,1}$, such that the complement of their
union in each component of  $\Sigma^{\infty}_{d,r} \setminus \Sigma_{d,r,1}$ has $2$ connected components.
\end{itemize}
These pieces must be {\em compatible} in the following sense:
first, the traces of the $d$-prerigid structure on each $d$-leg pants (i.e. the intersections with pants) are made up of $d+1$ connected components, called {\em seams}; secondly, each boundary component of the pants intersects with exactly two components of the seams at two distinct points; thirdly,  the seams cut each pants into two components.  Note that these conditions imply that each component is homeomorphic to a disk. One says then that the suited  $d$-pants decomposition and the $d$-prerigid structure are {\em subordinate} to the  $d$-rigid structure.
\item By construction, $\Sigma^{\infty}_{d,r}$ is naturally equipped with a suited $d$-pants decomposition, which will be referred to below as the {\em canonical suited $d$-pants decomposition}. We also fix a $d$-prerigid structure on $\Sigma^{\infty}_{d,r}$ (called the {\em canonical $d$-prerigid structure}) compatible with the canonical suited $d$-pants decomposition. See Figure \ref{fig-inf-surf-seam}. Using the puncture model,  the seams of the canonical $d$-prerigid structure are just the intersections of $[-1,1]\times \{0\}$ with each $d$-pants. The resulting $d$-rigid structure is called the {\em canonical $d$-rigid structure} on $\Sigma^{\infty}_{d,r}$. Very importantly, for each admissible subsurface, the canonical $d$-rigid structure induces an order on the admissible boundaries.  In  Figure \ref{fig-inf-surf-seam}, the induced order on the admissible loops are counterclockwise. Using the puncture model, the admissible loops are ordered from left to right.

\item The seams cut each component of $\Sigma^{\infty}_{d,r} \setminus \Sigma_{d,r,1}$ into two pieces, we choose the front piece in each component and these $r$ pieces together form  the {\em visible side} of $\Sigma^{\infty}_{d,r}$.
\item A suited $d$-pants decomposition (resp. $d$-(pre)rigid structure) is {\em asymptotically
    trivial} if outside a compact subsurface of $\Sigma^{\infty}_{d,r}$, it coincides with the
  canonical suited $d$-pants decomposition (resp. canonical $d$-(pre)rigid structure).
\end{enumerate}
\end{defn}

\begin{rem}\label{rem-comp-disk}
It is important that the seams cut each $d$-pants into two components and each component is homeomorphic to a disk as the mapping class group of a disk is trivial.
\end{rem}

\begin{defn}\label{defn-asym-rgd}
 Let $\Sigma^{\infty}_{d,r} $ and $ \bar{\Sigma}^{\infty}_{d,r'}$ be two surfaces with $d$-rigid structure and let $\varphi:  \Sigma^{\infty}_{d,r} \to \bar{\Sigma}^{\infty}_{d,r'}$ be a homeomorphism. One says that
$\varphi$ is {\em asymptotically rigid} if there
exists an admissible subsurface $A\subset \Sigma^{\infty}_{d,r}$ such
  that:
 \begin{enumerate}
 \item $\varphi(A)$ is also admissible in $\bar{\Sigma}^{\infty}_{d,r'}$,
 \item $\varphi \mid_A$ maps the based boundaries to based  boundaries, admissible loops to admissible loops and  
 \item the restriction of $\varphi: \Sigma^{\infty}_{d,r} \setminus A\to \bar{\Sigma}^{\infty}_{d,r'}\setminus\varphi(A)$
 is {\em rigid}, meaning that it respects the traces of the canonical $d$-rigid structure, mapping the suited $d$-pants
 decomposition into the suited $d$-pants decomposition, the seams into the
 seams, and the visible side into the visible side. 
 \end{enumerate}
 If we drop the condition that  $\varphi$ should map the visible side into the visible side, $\varphi$ is called {\em asymptotically quasi-rigid}. The surface $A$ is called a {\em support} for $\varphi$.  
\end{defn}

\begin{rem}
We are not using the word ``support" in the usual sense, as the map outside the support defined above might well not being the identity, but the map is uniquely determined up to isotopy by Remark \ref{rem-comp-disk}.
\end{rem}

\begin{rem}
In \cite[Definition 2.3]{FK04}, they do not actually require that the support must contain the base. This will not make a difference, as one can always enlarge the support so that it contains the base.
\end{rem}

\begin{rem} \label{rem-iso-rgs}
The surface $\Sigma^{\infty}_{d,r+d-1}$ can be identified with the surface $\Sigma^{\infty}_{d,r}$  such that $\Sigma_{d,r+d-1,m} = \Sigma_{d,r,m+1}$ for any $m\geq 1$, and the $d$-rigid structure of  $\Sigma^{\infty}_{d,r}$ coincides with $d$-rigid on $\Sigma^{\infty}_{d,r+d-1}$ outside $\Sigma_{d,r,2}$. In this way, $\Sigma^{\infty}_{d,r}$ is asymptotically rigid homeomorphic to $\Sigma^{\infty}_{d,r+d-1}$ through the identity map.
\end{rem}

\begin{rem}\label{induced-asy-rig-str}
Let $\Sigma'$ be a subsurface of $\Sigma^{\infty}_{d,r}$ such that there exist an admissible subsurface $A$ of $\Sigma^{\infty}_{d,r}$ satisfying:
\begin{enumerate}
     \item $A \cap \Sigma'$ is a compact surface,
     
    \item The boundaries of $\Sigma'$ are disjoint from the admissible boundary components of $A$.
    
    \item If an admissible boundary component $L$ of $A$ is contained in $\Sigma'$, then the punctured disk component of  $\Sigma^{\infty}_{d,r}$ cutting along $L$ is also contained in $\Sigma'$.
\end{enumerate}

Then $\Sigma'$ has a naturally induced $d$-rigid structure. In fact, we can take $A\cap \Sigma'$ to be the base surface and $d$-rigid structure can simply be inherited from $\Sigma^{\infty}_{d,r}$.
One, of course, can choose different $A$ here which may give different induced $d$-rigid structure, but it is unique up to asymptotically rigid homeomorphism. 
\end{rem}

\subsection{Asymptotic mapping class groups surjecting to Higman--Thompson groups} Given a (possibly noncompact) surface $\Sigma$, recall the mapping class group of $\Sigma$ is defined to be  the group of isotopy classes of orientation preserving homeomorphisms of $\Sigma$ that fixes $\partial \Sigma$ pointwise, i.e.
$$\Map(\Sigma) = \Map(\Sigma,\partial \Sigma) := \text{Homeo}^+(\Sigma, \partial \Sigma)/ \text{Homeo}_0(\Sigma, \partial \Sigma).$$

With this, we can now define the asymptotic mapping class group and the half-twist asymptotic mapping class group.

\begin{defn}\label{groupB}
The \emph{asymptotic mapping class group} $\CB V_{d,r}(\Sigma)$ (resp. the \emph{half-twist asymptotic mapping class group} $\CH V_{d,r}(\Sigma)$) is the subgroup of $\Map(\Sigma_{d,r}^\infty)$ consisting of isotopy classes of asymptotically
rigid (resp. quasi-rigid) self-homeomorphisms of  $\Sigma^{\infty}_{d,r}$. When $\Sigma$ is the disk, we sometimes simply denote the group by $\CB V_{d,r}$ (resp. $\CH V_{d,r}$). 
\end{defn}

\begin{defn}
Let $A$ be an admissible subsurface of $\Sigma_{d,r}^\infty$, and $\Map(A)$ be its mapping class group which fixes the each boundary component pointwise. Each inclusion  $A\subseteq A'$ of admissible surfaces induces an injective embedding $j_{A,A'} : \Map(A) \to \Map(A')$. The collection forms a direct system whose direct limit we call the \emph{compactly supported pure mapping class group}, denoted by $\PMap_c(\Sigma_{d,r}^\infty)$. The group $\PMap_c(\Sigma_{d,r}^\infty)$ is naturally a subgroup of $\CB V_{d,r}(\Sigma)$ and we denote the inclusion map by $j$.
\end{defn}

\begin{defn}\label{defn-forest}
Let ${\CF}_{d,r}$ be the forest with $r$ copies of a rooted $d$-ary tree and $\CT_{d,r}$ be the rooted tree obtained from ${\CF}_{d,r}$ by adding an extra vertex to ${\CF}_{d,r}$ and $r$ extra edges each connecting this vertex to a root of a tree in ${\CF}_{d,r}$.   There is a natural projection
$q:\Sigma^{\infty}_{d,r} \to {\CT_{d,r}}$, such that the pullback of the root is $\Sigma_{d,r,0}$ and the pull back of the midpoints of any edges are admissible loops.
\end{defn}

Now any element in $\CB V_{d,r}(\Sigma)$ can be represented by an asymptotically rigid homeomorphism $\varphi:  \Sigma^{\infty}_{d,r} \to {\Sigma}^{\infty}_{d,r}$. In particular we have an admissible subsurface $A$ of $\Sigma^{\infty}_{d,r}$ such that $\varphi|_A :(A,\partial_b A)\to (\varphi(A),\varphi(\partial_b A))$ is a homeomorphism. Let $F_-$ be the smallest subforest of $\CF_{d,r}$ which contains $q(A) \cap \CF_{d,r}$, and $F_+$ be the smallest subforest of $\CF_{d,r}$ which contains $q(\varphi( A)) \cap \CF_{d,r}$. Note that $F_-$ and $F_+$ have the same number of leaves and their leaves are in one-to-one correspondence with the admissible loops of $A$ and $\varphi(A)$.  Now let $\rho$ be the map from leaves of $F_-$ to $F_+$ induced by $\varphi$. Together this defines an element $[(F_-,\rho,F_+)]\in V_{d,r}$. We call this map $\pi$. One can show $\pi$ is well defined. Similarly to \cite[Proposition 2.4]{FK04} and \cite[Proposition 4.2, 4.6]{AF17}, we now have the following proposition.

\begin{prop}\label{thm-ses-asm}
We have the short exact sequences:
$$ 1 \to \PMap_c(\Sigma_{d,r}^\infty) \xrightarrow{j} \CB V_{d,r}(\Sigma) \xrightarrow{\pi} V_{d,r} \to 1;$$
$$ 1 \to \PMap_c(\Sigma_{d,r}^\infty) \xrightarrow{j} \CH V_{d,r}(\Sigma) \xrightarrow{\pi} V_{d,r}(\BZ/2\BZ) \to 1.$$
\end{prop}

\begin{rem}
Here, as in \cite{AF17}, $V_{d,r}(\BZ/2\BZ) $ is the twisted version of the Higman--Thompson group where one allows flipping the subtree below every leaf. See for example \cite{BDJ17} for more information.
\end{rem}

\Proof
We will prove the proposition for $\CB V_{d,r}(\Sigma)$. The other case is essentially the same. First we show the map $\pi$ is surjective. Given any element  $[(F_-,\rho,F_+)] \in V_{d,r}$, let $T_-$ (resp. $T_+$) be the tree obtained from $F_-$  (resp. $F_+$)  by adding a single root on the top and $r$ edges connecting to each root of the trees in $F_-$ (resp. $F_+$). Furthermore, let  $T'_-$ (resp. $T'_+$) be the tree obtained from $F_-$  (resp. $F_+$)  by throwing away the leaves and the open half edge connecting to the leaves. Then let $A_- =q^{-1}(T'_+)$ and $A_+ =q^{-1}(T'_+)$. We have $A_-$ and $A_+$ are both admissible subsurfaces of $\Sigma_{d,r}^{\infty}$.  Now one can produce a homeomorphism $\varphi_0: A_- \to A_+$ which is identity on the based boundary and maps the admissible loops of $A_-$ to the admissible loops of $A_+$ following the information from $\rho$, mapping the visible part to the visible part for each admissible loop. From here, we extend $\varphi_0$ to a map $\varphi: \Sigma_{d,r}^\infty \to \Sigma_{d,r}^\infty$ such that $\varphi$ is a asymptotically rigid homeomorphism. 

If an element $g\in \CB V_{d,r}(\Sigma) $ is mapped to a trivial element $\pi(g) = [(F_-,\rho,F_+)] \in V_{d,r}$, then the two forests $F_-$ and $F_+$ are the same and the induced map $\rho$ is trivial. This means we can assume the support $A$ for the asymptotically rigid homeomorphism $\varphi_g$ corresponding to $g$ is the same as $\varphi(A)$ and $\varphi$ induces identity map on the admissible boundary components.  Thus $g \in \PMap_c(\Sigma_{d,r}^\infty)$. Finally, given any element $g\in \PMap_c(\Sigma_{d,r}^\infty)$, it is clear that $\pi \circ j(g) = 1$.
\qed
 
The mapping class group $\Map(\Sigma_{d,r}^\infty)$ has a natural 
quotient topology coming from the compact-open topology on $\text{Homeo}^+(\Sigma_{d,r}^\infty, \partial \Sigma_{d,r}^\infty)$. See \cite[Section 2.3, 4.1]{AV20} for more information. In \cite[Theorem 1.3]{AF17}, Aramayona and Funar showed that when $\Sigma$ is a closed surface, $\CH V_{2,1}(\Sigma)$ is dense in $\Map(\Sigma_{2,1}^\infty)$. We improve their result to the following.

\begin{thm}\label{thm-dense-asy}
The groups $\CB V_{d,r}(\Sigma)$ and $\CH V_{d,r}(\Sigma)$ are dense in the mapping class group $\Map(\Sigma_{d,r}^\infty)$.
\end{thm}
\Proof 
The proof in \cite[Section 6]{AF17} adapts directly to show that $\CH V_{d,r}(\Sigma)$ is dense in $\Map(\Sigma_{d,r}^\infty)$ and so we will not repeat it here. To show $\CB V_{d,r}(\Sigma)$ is also dense in $\Map(\Sigma_{d,r}^\infty)$, it suffices to show any element in $\CH V_{d,r}(\Sigma)$ can be approximated by a sequence of elements in $\CB V_{d,r}(\Sigma)$. Note first that any half Dehn twists around  admissible loops in $\Sigma_{d,r}^\infty$ lies in $\CH V_{d,r}(\Sigma)$. In fact, given an admissible loop $\alpha$, we can choose an admissible subsurface $A$ such that $\alpha$ is an admissible loop of $A$. Then the half Dehn twist around $\alpha$ is asymptotic quasi-rigid with support $A$, in fact,  it is identity on all the components of $\Sigma_{d,r}^\infty \setminus A$ except at the component containing $\alpha$  where it rotates 180 degree. Now given an asymptotic quasi-rigid homeomorphism $f$ of $\Sigma_{d,r}^\infty$ with support $A'$,  we can first compose $f$ with  half  Dehn twists around  those admissible loops of $A$ where $f$ restricted to the component below them switches the front and back. The composition now is an asymptotic rigid homeomorphism. Thus $\CH V_{d,r}(\Sigma)$ can be generated by $\CB V_{d,r}(\Sigma)$ and half Dehn twists around the admissible loops in  $\Sigma_{d,r}^\infty$, it suffices now to show that any half Dehn twists around admissible loops in $\Sigma_{d,r}^\infty$ can be approximated by a sequence of elements in $\CB V_{d,r}(\Sigma)$. Given an admissible loop $L$, let $h_L$ be a half Dehn twist at $L$. We will construct a sequence of elements $x_i\in \CB V_{d,r}(\Sigma)$ such that for any compact subset $K$ of $\Sigma_{d,r}^\infty$, there exists $N$ such that for any $j\geq N$, $x_j$ and $h_L$ coincide on $K$. Recall we have the map $q: \Sigma_{d,r}^\infty \to \CT_{d,r}$ (cf. Definition \ref{defn-forest}) such that the admissible loops are mapped to edge middle points in $\CT_{d,r}$. Now consider those admissible loops such that their image under $q$ lying below $q(L)$ have distance $i$ to $q(L)$. Note that there are $d^i$ such admissible loops. We list them as $L_{i,1},\cdots, L_{i,d^i}$. Let $h_{L_{i,k}}$ be the half Dehn twists around $L_{i,k}$ and let $x_i =h_L h_{L_{i,1}} \cdots h_{L_{i,d^i}}$, then $x_i\in \CB V_{d,r}(\Sigma)$ and the sequence $\{x_i\}$ has the desired property. 
\qed

Now recall by Remark~\ref{rem-emb-model-gen} that $\Sigma_{d,r}^\infty$ is homeomorphic to $\Sigma \setminus \CC$ for any $d$ and $r$, hence we have our first theorem stated in the introduction.

%\begin{cor}\label{thm-dense}
%Let $\Sigma$ be any compact surface and $\CC$ be a Cantor set which lies in the interior of a disk in $\Sigma$. Then the mapping class group $\Map(\Sigma\setminus \CC)$ contains the following two families of dense subgroups: the asymptotic mapping class groups $\CB V_{d,r}(\Sigma)$ which surject to the Higman--Thompson group $V_{d,r}$, and  the half-twist asymptotic mapping class groups $\CH V_{d,r}(\Sigma)$  which surject to the symmetric Higman--Thompson group $V_{d,r}(\BZ/2\BZ)$.    
%\end{cor}

\subsection{The asymptotic mapping class group of the disk punctured by the Cantor set}
In the last subsection, we want to identify the asymptotic mapping class group $\CB V_{d,r}(D)$ with the oriented ribbon Higman--Thompson groups $RV^+_{d,r}$ and the half-twist asymptotic mapping class group $\CH V_{d,r}(D)$ with the ribbon Higman--Thompson group $RV_{d,r}$. The following lemma appears in \cite[Section 2]{BT12} without a proof, so we provide the details here. Note that what they call the pure ribbon braid group is the oriented pure ribbon braid group in our definition (see Definition \ref{defn-rb}). 

\begin{lem}\label{lem-pbraid-ribb}
Let $D_k$ be the $(k+1)$-holed sphere. Then $\Map(D_k)$ can be naturally identified with the  pure oriented ribbon braid group $PRB_k^+$. 
\end{lem}
\Proof Note that $D_k$ can be identified with a disk with $k$ holes. Let $\partial_b$ denote the boundary of the disk. Let $\bar{D}_k$ be a disk with $k$ punctures obtained from $D_k$ by attaching one punctured disk to each hole. The induced map $\CC ap:\Map(D_k)\to \PMap(\bar{D}_k)$ is the capping homomorphism. Note that $\PMap(\bar{D}_k) \cong PB_k$. Now applying \cite[Proposition 3.19]{FaMa11} and the fact that the Dehn twists around the holes of $D_k$ commute, one sees that the kernel $K$ is a free abelian group of rank $k$ generated by these $k$ Dehn twists.  Here the capping homomorphism splits. To prove this, we first embed $PB_k$ into $PRB_k^+$ by viewing the pure braid group of $k$-strings as the set of ribbon braids on $k$ bands such that the bands have no twists. We can think of $D_k$ as being embedded into $\mathbb{R}^2$ with $\partial_b$ as the unit circle and the $k$ holes in $D_k$ equally distributed inside $\partial_b$ along the $x$-axis. The intersections of these holes with the $x$-axis gives $k$ sub-intervals of the $x$-axis denoted $I_1,\cdots,I_k$.  We now put the bands representing a pure braid $x\in PB_k\leq PRB_k^+$ in $D\times [0,1]$ which starts and ends at $I_1,\cdots,I_k$. Note that the bands here will not twist at all. Now we comb the bands straight from bottom to top.  This induces a homeomorphism of $D_k\times \{0\}$ and hence an element in the mapping class group  $\Map(D_k)$. One checks that this map is a group homomorphism and injective. Since $PB_k$ acts on $K$ trivially, we have $\Map(D_k) \cong K \times PB_k \cong \BZ^k\times PB_k \cong PRB_k^+$ where the number of Dehn twists around each boundary component is naturally identified with the number of full twists on each bands.
\qed

To promote Lemma \ref{lem-pbraid-ribb} such that it works for the ribbon braid group, we need some extra terminology. As in the proof of  Lemma \ref{lem-pbraid-ribb}, we identify $D_k$ with the unit disk in $\mathbb{R}^2$ with $k$ small disks whose centers are equally distributed on the $x$-axis removed. The $x$-axis cuts the boundary loops of each deleted disk into two components, providing a cell structure on the loops. We will call the part that lies above the $x$-axis the \emph{visible part}. We define the \emph{rigid mapping class group} $\RMap_+(D_k)$ of $D_k $ to be the isotopy classes of homeomorphisms of $D_k$ which fix $\partial_b D_k$ pointwise and map the visible part of the holes to the visible part of the holes. Note  elements in $\partial_b D_k$ are allowed to map one boundary hole to another just as in the definition of the asymptotic mapping class group. If we only assume the cell structure on the loops has to be preserved, the resulting group is called \emph{quasi-rigid mapping class group} $D_k$ and denoted by $\RMap(D_k)$. With these preparations, the following lemma is now clear.

\begin{lem}\label{lem-braid-ribb}
There is a natural isomorphism between the oriented ribbon braid group $RB_k^+$ and $\RMap_+(D_k)$ (resp. between the ribbon braid group $RB_k$ and $\RMap(D_k)$).
\end{lem}
\Proof 
As in the proof of Lemma \ref{lem-pbraid-ribb}, we put the element in the (oriented) ribbon braid group between $D\times [0,1]$, then we comb the bands straight from bottom to top which gives the corresponding element in $\RMap_+(D_k)$ (resp. $\RMap(D_k)$ ).
\qed

Given two admissible subsurfaces $A$ and $A'$  of $D^\infty_{d,r}$ (possibly with different $r$) with $k$ admissible boundary components, we want to fix a canonical way to identify a homeomorphism $f:A\to A'$ as an element in the ribbon braid group. Note that each boundary loop except the base one inherits a visible side from $D^\infty_{d,r}$. We will use the puncture model for  $D^\infty_{d,r}$ going forward.

As above, let $D_k$ be the subsurface of $D$ which is the compliment of $k$ disjoint open disks with centers at $a_i = \frac{2i-k-1}{k+1} $ of radius $2^{-k}$ for $1\leq i\leq k$. Now given any admissible subsurface $A_k$ of $D^\infty_{d,r}$ with $k$ many admissible boundaries, denote the centers from left to right by $c_i \in [0,1]\times \{0\}$, $1\leq i\leq k$ with radius $r_1,r_2\cdots, r_{k}$. Now we define an isotopy $\CN_{A_k}:D\times [0,1] \to D$ such that $\CN_{A_k,0} =\id_D$ and $\CN_{A_k,1}$ maps $A_k$ to $D_k$ via a homeomorphism. We first shrink the admissible boundary loops of $A_k$ such that they have radius $r$, where $r = \min \{r_1,\cdots,r_k,2^{-k}\}$. Then we  isotope $A_k$ by moving the centers  $c_i$ to $a_i$ along $ [0,1] \times \{0\}$ in $D$. And in the last step we enlarge the radius one by one to $2^{-k}$. The following lemma is now immediate.

\begin{lem}\label{lem-home-brd}
Let  $\phi:D_{d,r}^\infty \to D_{d,r}^\infty$ be an asymptotically rigid (resp. quasi-rigid) homeomorphism which is supported on the admissible subsurface $A_k$. Denote $A_k' =\phi(A_k)$, then
\begin{enumerate}
    \item $\mathfrak{r}_\phi = \CN_{A_k',1} \circ \phi|_{A_k} \circ \CN_{A_k,1}^{-1}: D_k\to D_k$ gives an element in the oriented ribbon braid group $RB^+_k$ (resp. the ribbon braid group $RB_k$). Conversely, given an element $\mathfrak{r}\in RB^+_k$ (resp. $RB_k$), we have an asymptotic rigid (resp. quasi-rigid) homeomorphism which is unique up to isotopy,  supported on $A_k$,  and map $A_k$ to $A_k'$.
    
    \item let $A_{k+d}$ be the admissible subsurface of $D_{d,r}^\infty$ obtained from $A_k$ by adding a $d$-leg pants and $ \phi(A_{k+d})= A_{k+d}'$. Then the associated  oriented ribbon braid  (resp. the ribbon braid) of $\phi$ can be obtained from $\mathfrak{r}_{\phi}$ by splitting the corresponding band into $d$ bands. Conversely, if we split one band of the ribbon braids into $d$ bands, the isotopy class of the corresponding asymptotic  rigid (resp. quasi-rigid) homeomorphism does not change.  
\end{enumerate}
\end{lem}

Note that for any $d\geq2,r\geq1$, we have an natural embedding $\iota_{{d,r}}: D_{d,r}^\infty \to D_{d,r+1}^\infty$ by mapping the rooted boundaries of $D_{d,r}^\infty $ to the first $r$ rooted boundaries according to the order. This induces an embedding of groups $i_{\CH,{d,r}}: \CH V_{d,r} \to \CH V_{d,r+1}$, $i_{\CB,{d,r}}: \CB V_{d,r} \to \CB V_{d,r+1}$.
On the other hand, we also have natural embeddings  $i_{R,d,r}: RV_{d,r} \to RV_{d,r+1} $ and $i_{ R^+,d,r} :  RV^+_{d,r} \to  RV^+_{d,r+1}$ induced by inclusion of roots. We have the following.

\begin{thm}\label{thm-iden-asym-ribb}
     There exist isomorphisms $f_{d,r}:\CH V_{d,r} \rightarrow RV_{d,r}$ such that $f_{d,r+1}  \circ i_{\CH,d,r}  = i_{\CH,d , r+1} \circ f_{d,r}$. Restricting to the subgroups $\CB V_{d,r}$, one gets isomorphisms $f_{d,r}:\CB V_{d,r} \rightarrow RV^+_{d,r}$ with the same property.
\end{thm}

\Proof
Since the two cases are parallel, we will only prove the theorem for $\CB V_{d,r}$.  We will define two maps $f_{d,r}:\CB V_{d,r} \rightarrow RV^+_{d,r} $ and $g_{d,r}:RV^+_{d,r} \to \CB V_{d,r}$ such that $f_{d,r} \circ g_{d,r} =\id$ and $g_{d,r}\circ f_{d,r} = \id$.

Given an element $x \in \CB V_{d,r}$, we can define $f_{d,r}$ as follows. Let $\varphi_x$ be an asymptotically rigid homeomorphism of $D_{d,r}^{\infty}$ representing $x$ with support $A_k$, where $k$ is the number of admissible loops.  By Proposition \ref{thm-ses-asm}, $\pi(x)$ provides an element $[F_-,\sigma,F_+]$ in the Higman--Thompson group $V_{d,r}$, where $F_-$ and $F_+$ are $(d,r)$-forests with $k$ leaves.  But what we want is a ribbon braid connecting the $k$ leaves. For this we simply apply Lemma \ref{lem-home-brd} (1) to the map $\varphi_x$ with support $A_k$, denote the corresponding element in $RB^+_k$ by $\mathfrak{r}_{\varphi_x}$. We define $f_{d,r}(x) = [F_-,\mathfrak{r}_{\varphi_x},F_+]$.

Now given $y\in RV^+_{d,r}$, one can define an element in $ \CB V_{d,r}$ as follows. Suppose $(F_-,\mathfrak{r}_y, F_+)$ is a representative for $y$, where $F_-$ and $F_+$ are $(d,r)$-forests  and $\mathfrak{r}$ is a ribbon braid between the leaves of $F_-$ and $F_+$. Add a root to $F_-$ (resp. $F_+$) with an edge connecting to the root of each tree in $F_-$ (resp. $F_+$) and then throw away the open half edge connecting to the leaves. Denote the resulting tree by $T_-$ (resp. $T_+$).  Now $q^{-1}(T_-)$, $q^{-1}(T_+)$ gives us two admissible subsurfaces $A_k$, $A_k'$ in $D^{\infty}_{d,r}$  where $k$ is the number of leaves for $F_-$.  And by Lemma \ref{lem-home-brd} (1), the ribbon element $\mathfrak{r}_y$ in $RB^+_k$ give us an asymptotic rigid  homeomorphism $\psi_y$ with support $A_k$ and maps $A_k$ to $A_k'$. 

Now one can check that $f_{d,r} \circ g_{d,r} =\id$ and $g_{d,r}\circ f_{d,r} = \id$. Therefore, the two groups are isomorphic. The fact that the diagram commutes is immediate from the definition.
\qed

\section{Homological stablity of ribbon Higman--Thompson groups}\label{sec:homstab}
In this section, we show the homological stability for oriented ribbon Higman--Thompson groups and explain at the end how the same proof applies to the ribbon Higman--Thompson groups. 

%and then indicating at the end how the proof can be adapted to show the homological stability for the braided  Thompson groups.

\subsection{Homogeneous categories and homological stability}

In this subsection, we review the basics of homogeneous categories and refer the reader to \cite{RWW17} for more details. Note that we adopt their convention of identifying objects of a category with their identity morphisms.

\begin{defn}[{\cite[Definition 1.3]{RWW17}}]
A monoidal category $(\mathcal{C}, \oplus, 0)$ is called \emph{homogeneous} if $0$ is initial in $\mathcal{C}$ and if the following two properties hold.
\begin{inlinecond}{\textbf{H1}} $\Hom(A,B)$ is a transitive $\Aut(B)$-set under postcomposition.\end{inlinecond}
\begin{inlinecond}{\textbf{H2}} The map $\Aut(A) \rightarrow \Aut(A \oplus B)$ taking $f$ to $f \oplus \id_B$ is injective with image
\[\Fix(B) := \{ \phi\in \Aut(A\oplus B) \mid \phi\circ (\imath_A\oplus \id_B) = \imath_A \oplus \id_B\}\]
where $\imath_A\colon 0 \to A$ is the unique map.\end{inlinecond}
\end{defn}

\begin{defn}[{\cite[Definition 1.5]{RWW17}}]
Let $(\mathcal{C}, \oplus, 0)$ be a monoidal category with $0$ initial. We say that $\CC$ is \emph{prebraided} if its underlying groupoid is braided and for each pair of objects $A$ and $B$ in $\CC$, the groupoid braiding $b_{A,B} :A \oplus B \to B\oplus A$ satisfies 
$$ b_{A,B} \circ (A\oplus \imath_B) = \imath_B\oplus A : A\to B\oplus A.$$
\end{defn}

\begin{defn} \cite[Definition 2.1]{RWW17}
Let $(\CC,\oplus,0)$ be a monoidal category with $0$ initial  and $(A,X)$ a pair of objects in $\CC$. Define $W_n(A,X)_\bullet$ to be the semi-simplicial set with set of $p$-simplices
\[ W_n(A,X)_p := \Hom_\CC(X^{\oplus p+1}, A\oplus X^{\oplus n})\]

and with face map
$$d_i : \Hom_\CC(X^{\oplus p+1}, A \oplus X^{\oplus n}) \to \Hom_\CC (X^{\oplus p}, A \oplus X^{\oplus n}) $$
defined by precomposing with $X^{\oplus i} \oplus \imath_X \oplus X^{\oplus p-i}$

Also call the category $\mathcal{C}$ satifies $LH_3$ at a pair of objects $(A,X)$ with \emph{slope} $k\ge2$ if:  
\begin{inlinecond}{\textbf{LH3}}\label{cond:H3}
\textit{For all $n \geq 1$, $W_n(A,X)_\bullet$ is $(\frac{n-2}{k})$-connected.}\end{inlinecond}
\end{defn}

Quite often, we can reduce the semi-simplicial complex to a simplicial complex. 

\begin{defn}[{\cite[Definition 2.8]{RWW17}}] \label{defn:Sn(X,A)}
Let $A,X$ be objects of a homogeneous category $(\mathcal{C}, \oplus,0)$. For $n \geq 1$, let $S_n(A,X)$ denote the simplicial complex whose vertices are the maps $f\colon X \to A\oplus X^{\oplus n}$ and whose $p$-simplices are $(p+1)$-sets $\{ f_0, \dots, f_p\}$ such that there exists a morphism $f\colon X^{\oplus p+1} \to A\oplus X^{\oplus n}$ with $f \circ i_j = f_j$ for some order on the set, where
\[ i_j =  \imath_{X^{\oplus j}} \oplus \id_X \oplus \imath_{X^{\oplus p-j}}  \colon X = 0 \oplus X \oplus 0 \longrightarrow X^{\oplus p+1}. \]
\end{defn}

\begin{defn}
Let $Aut(A\oplus X^{\oplus \infty})$ be the colimit of 
$$ \cdots \xrightarrow[]{-\oplus X} Aut(A\oplus X^{\oplus n}) \xrightarrow[]{-\oplus X} Aut(A\oplus X^{\oplus n+1}) \xrightarrow[]{-\oplus X} Aut(A\oplus X^{\oplus n+2})\xrightarrow[]{-\oplus X}\cdots $$
Then any $Aut(A\oplus X^{\oplus \infty})$-module $M$ may be considered as an $Aut(A\oplus X^{\oplus n})$-module for any $n$, by restriction, which we continue to call $M$. We say that the module $M$ is abelian if the action of $Aut(A\oplus X^{\oplus \infty})$ on $M$ factors through the abelianizations of $Aut(A\oplus X^{\oplus \infty})$, or in other words if the derived subgroup of $Aut(A\oplus X^{\oplus \infty})$ acts trivially on $M$.  
\end{defn}

We are now ready to quote the theorem that we will use.
\begin{thm}[{\cite[Theorem 3.4]{RWW17}}] \label{thm:hom stab}
Let $(\mathcal{C}, \oplus,0)$ be a pre-braided homogeneous category satisfying \ref{cond:H3} for a pair $(A,X)$ with slope $k \geq 3$. Then for any abelian $Aut(A\oplus X^{\oplus \infty})$-module $M$ the map
\[H_i(\Aut(A\oplus X^{\oplus n}); M )  \longrightarrow H_i(\Aut(A\oplus X^{\oplus n+1}); M)\]
induced by the natural inclusion map is surjective if  $i \leq \frac{n-k+2}{k}$, and injective if  $i \leq \frac{n-k}{k}$.
\end{thm}

\subsection{Homogeneous category for the groups $RV^+_{d,r}$}

The purpose of this section is to produce a homogeneous category for proving homological stability of the ribbon Higman--Thompson groups $RV^+_{d,r}$.  Note that by Theorem \ref{thm-iden-asym-ribb}, it is the same as proving the asymptotic mapping class groups $ \CB V_{d,r}$  have homological stability. This allows us to define our homogeneous category geometrically. The category is similar to the ones produced in \cite[Section 5.6]{RWW17}. Essentially, we replace the annulus or M\"obius band by the infinite surface $D^{\infty}_{d,1}$.

\begin{figure}
\centering
\includegraphics[width=1\textwidth]{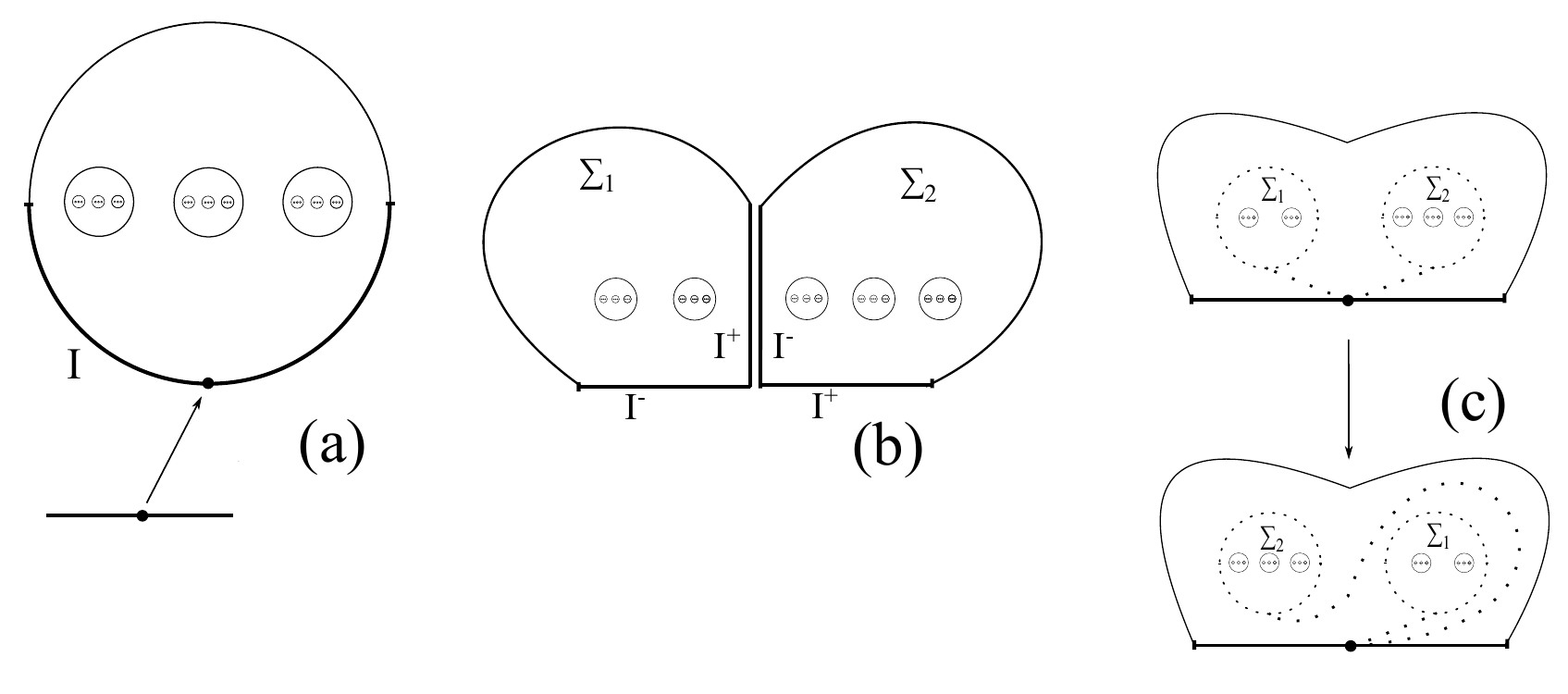}
\caption{The braided monoidal structure for the category $\CG_d$}
\label{fig-braid-monoid}
\end{figure}

Recall $D^{\infty}_{d,r}$ is an infinite surface equipped with a canonical asymptotic rigid structure with boundary component denoted $\partial_b D^{\infty}_{d,r}$. Let $I=[-1,1] \subset \partial_b D^{\infty}_{d,r}$ be an embedded interval as in Figure~\ref{fig-braid-monoid}(a). Let $I^- = [-1,0]$ and $I^+ = [0,1]$  be subintervals of $I$. Let $ D^{\infty}_{d,1} \oplus D^{\infty}_{d,1}$ be the boundary sum of two copies of $D^{\infty}_{d,1}$ obtained by identifying $I^+$ of the first copy with $I^-$ of the second copy. Inductively, we could define similarly $\oplus_r D^{\infty}_{d,1}$ for any $r\geq 0$.  Here $\oplus_0 D^{\infty}_{d,1}$ is just the standard disk $D$. Abusing notation, when referring to $I^-$ and $I^+$ on $\oplus_r D^{\infty}_{d,1}$, we will mean the two copies of $I^-$ and $I^+$ which remain on the boundary.  Thus we have an operation $\oplus$ on the set $\oplus_r D^{\infty}_{d,1}$ for any $r\geq 0$. See Figure \ref{fig-braid-monoid}(b) for a picture of $(\oplus_2 D^{\infty}_{d,1}) \oplus (\oplus_3 D^{\infty}_{d,1})$.  In fact, we have $ ((\oplus_r D^{\infty}_{d,1}), \oplus) $ is the free monoid generated by $D^{\infty}_{d,1}$.    Note that $\oplus_r D^{\infty}_{d,1}$ has a naturally induced $d$-rigid structure and we can identify it with $D^{\infty}_{d,r}$, which will be of use to us later.

We can now define the category $\CG_d$ to be the monoidal category with objects $\oplus_r D^{\infty}_{d,1}$, $r\ge 0$, $\oplus$ as the operation, and $D$ as the $0$ object. So far it is the same as defining the objects as the  natural numbers and addition as the operation. When $r=s$, we define the morphisms $\Hom (\oplus_r D^{\infty}_{d,1}, \oplus_s D^{\infty}_{d,1}) = \CB V_{d,r}$  which is the group of isotopy classes of  asymptotically rigid homeomorphisms of $ D^{\infty}_{d,r}$; when $r\neq s$, let $\Hom (\oplus_r D^{\infty}_{d,1}, \oplus_s D^{\infty}_{d,1}) =\emptyset $. Note that we did not universally define the morphisms to be the sets of isotopy classes of asymptotically rigid homeomorphisms as we want our category to satisfy cancellation, i.e., if $A\oplus C = A$ then $C = 0$, see \cite[Remark 1.11]{RWW17} for more information. The category $\CG_d$ has a natural braiding as in the usual braid group case, see Figure \ref{fig-braid-monoid}(c). 

Now applying \cite[Theorem 1.10]{RWW17}, we have a pre-braided homogeneous category $U \CG_d$. The category $U \CG_d$ has the same objects as $\CG_d$ and morphisms defined as following: For any $s \leq r$, a morphism in $\Hom (\oplus_s D^{\infty}_{d,1}, \oplus_r D^{\infty}_{d,1})$ is an equivalence class of pairs $(\oplus_{r-s} D^{\infty}_{d,1},f)$ where  $f: (\oplus_{r-s} D^{\infty}_{d,1}) \oplus (\oplus_{s} D^{\infty}_{d,1} )\to \oplus_r D^{\infty}_{d,1}$ is a morphism in $\CG_d$ and  $(\oplus_{r-s} D^{\infty}_{d,1},f) \sim (\oplus_{r-s} D^{\infty}_{d,1},f')$ if there exists an isomorphism $g: \oplus_{r-s} D^{\infty}_{d,1}\to \oplus_{r-s} D^{\infty}_{d,1} \in \CG_d$ making the diagram commute up to isotopy.

\[\begin{tikzcd}
(\oplus_{r-s} D^{\infty}_{d,1}) \oplus (\oplus_s D^{\infty}_{d,1}) \arrow[d,"g~\oplus~ \id_{\oplus_s D^{\infty}_{d,1}}"'] \arrow[r, "f"] &  \oplus_r D^{\infty}_{d,1} \\
(\oplus_{r-s} D^{\infty}_{d,1}) \oplus (\oplus_s D^{\infty}_{d,1}) \arrow[ur, "f'"']  &
\end{tikzcd}
\]
 We write $[\oplus_{r-s} D^{\infty}_{d,1},f]$ for such an equivalence class.  Now by Theorem \ref{thm:hom stab}, to prove homological stability for the oriented ribbon Higman--Thompson groups, we  only need to verify that the category $\mathcal{G}_d$ satisfies Condition \ref{cond:H3} at the pair $(D,D^{\infty}_{d,1})$. In fact, by Theorem \ref{thm-iden-asym-ribb}, proving the oriented ribbon Higman--Thompson groups satisfy homological stability is the same as proving that the the asymptotic mapping class groups $ \CB V_{d,r}$ satisfy homological stability. Now consider the family of groups 
 \[\Aut(A\oplus X) \hookrightarrow \Aut(A\oplus X^{\oplus 2}) \hookrightarrow \Aut(A\oplus X^{\oplus 2})\hookrightarrow \cdots \hookrightarrow \Aut(A\oplus X^{\oplus n}) \hookrightarrow \cdots \]
 Where $A = D$ and $X= D^{\infty}_{d,1}$. By definition, this gives rise to the family of groups $\CB V_{d,1} \hookrightarrow \CB V_{d,2} \hookrightarrow \cdots \hookrightarrow \CB V_{d,n} \hookrightarrow \cdots$. Now we have shown that the category $(\mathcal{G}_d, \oplus, D)$ is a prebraided homogenuous category, so  by Theorem \ref{thm:hom stab}, it suffices to verify Condition \ref{cond:H3} at the pair $(D,D^{\infty}_{d,1})$ to prove our homological stablity result. As a matter of fact, we will show that $W_r(D, D^{\infty}_{d,1})_\bullet$ is $(r-3)$-connected in the next subsection. First, let us further characterize the morphisms in $U\CG_d$.  Call $0 = I^- \cap I^+$ the basepoint of $\oplus_r D^{\infty}_{d,1}$.

\begin{defn}\label{defn-qr-emb}
Given $s < r$, an injective map $\varphi: (\oplus_s D^{\infty}_{d,1},I^+)\to  (\oplus_r D^{\infty}_{d,1},I^+)$ is called an \emph{ asymptotically rigid embedding} if it satisfies the following properties:
\begin{enumerate}
    
    \item $\varphi(\partial D^{\infty}_{d,s}) \cap \partial D^{\infty}_{d,r} = I^+$. 
     
        \item $\varphi$ maps $\oplus_s D^{\infty}_{d,1}$ homeomorphically to $\varphi(\oplus_s D^{\infty}_{d,1})$ and there exists an admissible surface $A\subset \oplus_s D^{\infty}_{d,1}$ such that  $\varphi: \oplus_s D^{\infty}_{d,1} \setminus A\to \varphi (\oplus_s D^{\infty}_{d,1}) \setminus \varphi(A)$  is rigid. 
    
    \item \label{defn-qr-emb-cplt} The closure of the complement of $\varphi(\oplus_s D^{\infty}_{d,1})$ in $\oplus_r D^{\infty}_{d,1}$ with its induced $d$-rigid structure is asymptotically rigidly homeomorphic to $\oplus_{r-s} D^{\infty}_{d,1}$. 
    
\end{enumerate}
\end{defn}

\begin{lem}
For $s<r$, the equivalence classes of pairs $[\oplus_{r-s} D^{\infty}_{d,1},f]$ are in one-to-one correspondence with the isotopy classes of asymptotically rigid embeddings of $(\oplus_s D^{\infty}_{d,1},I^+)$ into $(\oplus_r D^{\infty}_{d,1},I^+)$. 
\end{lem}
\begin{rem}
Here isotopies are carried out among asymptotically rigid embeddings.
\end{rem}

\Proof Given an equivalence class of a pair $[\oplus_{t-s} D^{\infty}_{d,1},f]$, we have the restriction map $f\mid_{\oplus_s D^{\infty}_{d,1}}$ is an asymptotically rigid embedding. Any two equivalence classes of pairs will induce the same map $f\mid_{\oplus_s D^{\infty}_{d,1}}$, hence we have a well-defined map from the set of equivalence pairs to the set of isotopy classes of asymptotically rigid embeddings.

We produce the inverse of the restriction map as follows. If we have an asymptotically rigid embedding $\varphi : (\oplus_s D^{\infty}_{d,1},I^+ )\to  (\oplus_r D^{\infty}_{d,1},I^+)$, by part \ref{defn-qr-emb-cplt} of Definition \ref{defn-qr-emb}, we also have an asymptotically rigid homeomorphism $\phi: C \to \oplus_{r-s} D^{\infty}_{d,1}$ where $C$ is the closure of the complement of $\varphi(\oplus_s D^{\infty}_{d,1})$ in $\oplus_{r} D^{\infty}_{d,1}$. Up to isotopy, we can assume $\phi^{-1}\mid_{I^+}$  coincides with $\varphi\mid_{I^-}$. Now define a map $\bar{f}: (\oplus_{r-s} D^{\infty}_{d,1}) \oplus (\oplus_s D^{\infty}_{d,1}) \to  \oplus_{r} D^{\infty}_{d,1}$ by $\bar{f}\mid_{\oplus_{r-s} D^{\infty}_{d,1}} =  \phi^{-1}$ and $\bar{f}\mid_{\oplus_s D^{\infty}_{d,1}} = \varphi$. One can check that $\bar{f}$ is an asymptotically rigid homeomorphism. Then $(\oplus_{r-s} D^{\infty}_{d,1},\bar f)$ gives a representative of an equivalence class of pairs. 
\qed

\subsection{Higher connectivity of the  complex $W_r(D,D^{\infty}_{d,1})_\bullet$}
We want to prove that the complex $W_r(D,D^{\infty}_{d,1})_\bullet$ is highly connected, see the diagram on Figure~\ref{fig:sumrelations}  and the paragraph preceeding it for an outline of our general strategy. 

\begin{rem}\label{rem:firstchange}
As explained in the proof of \cite[Lemma 5.21]{RWW17}, a simplex of $S_r(D,D^{\infty}_{d,1})$ has a canonical ordering on its vertices induced by the local orientation of the surfaces near the parameterized interval in their based boundary. Thus the geometric realization $|W_r(D,D^{\infty}_{d,1})_\bullet|$ is homeomorphic to $S_r(D,D^{\infty}_{d,1})$. 
\end{rem}

Our first step now is to simplify the complex $S_r(D,D^{\infty}_{d,1})$ further.

\begin{figure}[ht]
\centering
\includegraphics[width=.5\textwidth]{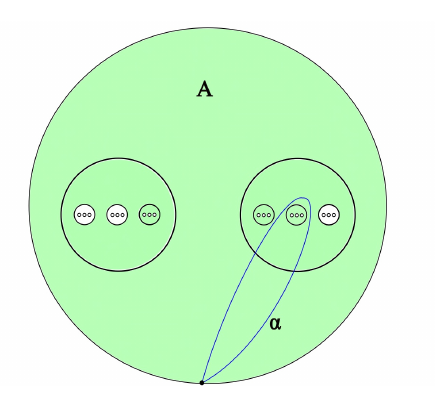}
\caption{The curve $\alpha$ is an asymptotically rigidly embedded loop, with the green shaded surface $A$ the corresponding admissible subsurface.}
\label{ad-loop}
\end{figure}

\begin{defn}
Given $r\geq 2$, we call a loop $\alpha: (I, \partial I)=([0,1],\{0,1\}) \to (\oplus_r D^{\infty}_{d,1},0)$ an \emph{asymptotically rigidly embedded loop} if there exists an asymptotically rigid embedding $\varphi: (D^{\infty}_{d,1},I^+ )\to  (\oplus_r D^{\infty}_{d,1},I^+)$ with $\varphi \mid_{(\partial D^{\infty}_{d,1},0)} =\alpha$ up to based isotopy. See Figure \ref{ad-loop} for a picture.
\end{defn}

\begin{rem}
When $r=1$, we just call a loop asymptotically rigidly embedded if it is isotopic to the boundary.
\end{rem}

\begin{lem}\label{lem-cha-asre}
When $r\geq 2$, a loop $\alpha: (I,\partial I) \to (\oplus_r D^{\infty}_{d,1},0)$ is isotopic to an asymptotically rigidly embedded loop if and only if there exists an admissible surface $A\subseteq \oplus_r D^{\infty}_{d,1}$ such that  the admissible loops of $A$ are disjoint from $\alpha$, the number of admissible loops of $A$ that lie in the disk bounded by $\alpha$ is $ 1 + a(d-1)$ for some $a\geq 0$ and there exist some admissible loops which do not lie inside the disk bounded by $\alpha$ up to isotopy.
\end{lem} 

\Proof It is clear that a loop which is isotopic to an asymptotically rigidly embedded loop has the properties given in the lemma. 

For the other direction, we can assume up to isotopy that $\alpha (I) \cap \partial (\oplus_r D^{\infty}_{d,1}) = I^+$. We know that $D^{\infty}_{d,r}$ is asymptotically rigidly homeomorphic to $D^{\infty}_{d,r+d-1}$, thus the surface bounded by the loop $\alpha$ is asymptotically rigidly homeomorphic to $D^{\infty}_{d,1}$. Therefore, the number of the boundary components bounded by the complement disk is $r-1\mod d-1$ and thus it is asymptotically rigidly homeomorphic to $D^{\infty}_{d,r-1}$. These two facts together imply $\alpha$ is an asymptotically rigidly embedded loop. 
\qed

Now we define the complex $U_r(D,D^{\infty}_{d,1})$ which is the surface version of the complex $U_r$ given in \cite[Section 2.4]{SW19}.

\begin{defn}
For $r \geq 1$, let $U_r(D,D^{\infty}_{d,1})$ denote the simplicial complex whose vertices are isotopy classes of asymptotically rigidly embedded loops and a set of vertices $\alpha_0,\cdots, \alpha_p$ forms a $p$-simplex if and only if any corresponding asymptotically rigid embeddings $\phi_0,\cdots,\phi_p$  form a $p$-simplex in $S_r(D,D^{\infty}_{d,1})$.
\end{defn}

We denote the canonical map from $S_r(D,D^{\infty}_{d,1})$ to $U_r(D,D^{\infty}_{d,1})$ by $\pi$. The following lemma follows directly from the definition. In fact, given a set of vertices $\alpha_0,\cdots, \alpha_p$, if they form a $p$-simplex, then  any corresponding asymptotically rigid embeddings $\phi_0,\cdots,\phi_p$  form a $p$-simplex in $S_r(D,D^{\infty}_{d,1})$. But this means for any $\psi_i \in \pi^{-1}(\alpha_i)$, $0\leq i\leq p$, the collection $\psi_0, \phi_1,\cdots, \psi_p$ forms a $p$-simplex.

\begin{lem}\label{lem-cpt-jn}
The map $\pi$ is a complete join.
\end{lem}

Now by Proposition \ref{prop-join-conn}, we only need to show that $U_r(D,D^{\infty}_{d,1})$ is highly connected.  Similar to \cite[Section 2.4]{SW19}, we will produce several other complexes closely related to $U_r(D,D^{\infty}_{d,1})$. We first have the following complex which is analogous to the complex $U_r^{\infty}$ in \cite[Definition 2.12]{SW19}.
 
\begin{defn}
  Let $U_r^{\infty}(D,D^{\infty}_{d,1})$ be the simplicial complex with vertices given by asymptotically rigidly embedded loops in  $\oplus_r D^{\infty}_{d,1}$ where $\alpha_0,\alpha_1,\cdots,\alpha_p$ form a $p$-simplex if the punctured disks bounded by them are pairwise disjoint (outside of the based point) and there exists  at least one admissible loop that does not lie in those disks.
\end{defn}

\begin{rem} \label{rem-Un}
\begin{enumerate}
    \item The $(r-2)$-skeleton of $U_r^{\infty}(D,D^{\infty}_{d,1})$ is the same as that of $U_r(D,D^{\infty}_{d,1})$. Notice though that $U_r^{\infty}(D,D^{\infty}_{d,1})$ is in fact infinite dimensional.
    
    \item Since $\oplus_rD^{\infty}_{d,1}$ is asymptotically rigidly homeomorphic to $\oplus_{r+d-1}D^{\infty}_{d,1}$, we have $U_r^{\infty}(D,D^{\infty}_{d,1})$ is isomorphic to $U_{r+d-1}^{\infty}(D,D^{\infty}_{d,1})$ as a simplicial complex. 
\end{enumerate}
\end{rem}

We also need the another complex which is the surface version of the complex $T_r^\infty$ given in \cite[Defintion 2.14]{SW19}.  For convenience, we will orient the admissible loops in $\oplus_r D^{\infty}_{d,1}$ such that they bound the punctured disk according to the orientation. 

\begin{defn}
An \emph{almost admissible loop} is a loop $\alpha: (I,\partial I) \to (\oplus_r D^{\infty}_{d,1},0)$ which is freely isotopic to one of the nonbased admissible loops.
\end{defn}

Note that by Lemma \ref{lem-cha-asre}, an almost admissible loop is an asymptotically rigidly embedded loop. 

\begin{defn}\label{defn-subcpxT}
Define the simplicial complex $T_r^{\infty}(D,D^{\infty}_{d,1})$ to be the full subcomplex of $U_r^{\infty}(D,D^{\infty}_{d,1})$ such that all its vertices are almost admissible loops. 
\end{defn}

Just as discussed in Remark \ref{rem-Un}, we have $T_r^{\infty}(D,D^{\infty}_{d,1})$ is in fact isomorphic to $T_{r+d-1}^{\infty}(D,D^{\infty}_{d,1})$ as a simplicial complex.

We now want to further characterise the almost admissible loops by building a connection to the usual arc complex. We let $A$ be the quotient $[0,2]/1\sim 2$. This corresponds to identifying the endpoint $1$ of the interval $[0,1]$ with the base point $1$ of the circle given by $[1,2]/1\sim 2$. 

\begin{defn}
 An injective continuous map  $L: (A,0) \to (D^{\infty}_{d,r},0)$ is called a \emph{lollipop} on the surface $D^{\infty}_{d,r}$ if $\alpha\mid_{[1,2]}$ is isotopic to an admissible loop in $D^{\infty}_{d,r}$ and $L\mid_{[0,1]}$ is an arc connecting the base point $0$ to the loop $L([1,2])$. The map $L\mid_{[0,1]}$ is called the \emph{arc part} of the lollipop $L$ and $L\mid_{[1,2]}$ is called the \emph{loop part}. See Figure \ref{jumping} where the blue curve is a lollipop.
\end{defn}

Lollipops are examples of what Hatcher-Vogtmann refer to as tethered curves \cite{HV17}.

\begin{lem} \label{lem:loop-arc}
The set of isotopy classes of almost admissible loops is in one-to-one correspondence with the set of isotopy classes of lollipops.  
\end{lem}

\Proof
We define a map $g$ from the isotopy classes of lollipops to the isotopy classes of almost admissible loops and show that the map is bijiective.

Given a lollipop $L: (A,0) \to (D^{\infty}_{d,r},0)$, we can map it to an almost admissible loop   $\alpha:[0,1] \to (D^{\infty}_{d,r},0) $ as follows. We define $\alpha(0)=0$ and let $\alpha(t)$ run parallel to $L$ outside the region bounded by $L$. The orientation of $\alpha$ is simply the one coincides with the loop part of $L$. Since $\alpha$ can be freely homotoped to the admissible loop $L\mid_{[1,2]}$, we have $\alpha$ is almost admissible. Any isotopy of $L$ induces an isotopy of $\alpha$, hence the map is well-defined.

Now we show $g$ is surjective. Given any almost admissible loop $\alpha:[0,1] \to (D^{\infty}_{d,r},0)$, let $A$ be the admissible loop which is freely isotopic to $\alpha$. Up to isotopy, we can assume that $A$ lies in the interior of the surface bounded by $\alpha$. Then the surface bounded by $\alpha$ and $A$ must be an annulus. From here one can produce an arc connecting the base point $0$ to a point in $A$. Together with $A$, this provides the lollipop. 

Finally, we argue that $g$ is injective. Suppose $L_1$ and $L_2$ are two lollipops such that $g(L_1)$ and $g(L_2)$ are isotopic, denote the isotopy by $f$. By the isotopy extension theorem (see for example \cite[Proposition 1.11]{FaMa11}) there exists an isotopy $F:D^{\infty}_{d,r} \times [0,1] \to D^{\infty}_{d,r}$ such that $F \mid_{D^{\infty}_{d,r} \times 0} = id_{D^{\infty}_{d,r}}$ and $F\mid_{g(L_1)\times [0,1]} =f$. In particular $F\mid_{D^{\infty}_{d,r} \times 1}$ maps the almost admissible loop $g(L_1)$ to the almost admissible loop $g(L_2)$. Hence $L_1$ is isotoped through $F$ to a lollipop which lies in a small neighborhood of $L_2$ and is bounded by the loop $g(L_2)$. Therefore, one can then isotope $L_1$ to $L_2$.
\qed

We now have the following definition of lollipop complex. 

\begin{defn}
The \emph{lollipop complex} $L_r^{\infty}(D,D^{\infty}_{d,1})$ has vertices as lollipops, and $p+1$ lollipops, $L_0,L_1,\cdots, L_p$, form a $p$-simplex if they are pairwise disjoint outside the base point $0$ and there exists at least one admissible loop which does not lie inside the disks bounded by the $L_i$s. 
\end{defn}

The following lemma is immediate from Lemma \ref{lem:loop-arc}. 

\begin{lem}\label{lem-id-T-L}
The complex $L_r^{\infty}(D,D^{\infty}_{d,1})$ is isomorphic to $T_r^{\infty}(D,D^{\infty}_{d,1})$ as a simplicial complex.
\end{lem}

\begin{lem}\label{lem-link-loli}
Given a $p$-simplex $\sigma$ in $L_r^{\infty}(D,D^{\infty}_{d,1})$, its link $\Lk({\sigma})$ is isomorphic to $L_{r_\sigma}^{\infty}(D,D^{\infty}_{d,1})$ for some $r_\sigma>0$.
\end{lem}
\Proof
By Lemma \ref{lem-id-T-L}, we can just prove the lemma for $T_r^{\infty}(D,D^{\infty}_{d,1})$. Let $\alpha_0,\alpha_1,\cdots,\alpha_p$ be the vertices of  $\sigma$ which are almost admissible loops. Up to isotopy, we can assume they are pairwise disjoint except at the basepoint $0$. Now let $C$ be the complement surface of  $\sigma$, whose based boundary is the concatenation of  $\alpha_p,\alpha_{p-1},\cdots,\alpha_0$ and $\partial D$. The surface $C$ has a naturally induced $d$-rigid structure. In particular, $C$ is asymptotically rigidly homeomorphic to $D^{\infty}_{d,r_\sigma}$ for some $r_\sigma>0$. Thus link $\Lk({\sigma})$ is isomorphic to $T_{r_\sigma}^{\infty}(D,D^{\infty}_{d,1})$.
\qed

\begin{figure}[ht]
 \[
   \begin{tikzpicture}[scale=0.85,every node/.style={thick,scale=0.85}]
    \matrix (m) [matrix of math nodes, column sep=1em, row sep=2em, text depth=.5em, text height=1em, ampersand replacement=\&]
    {
     \&  \&  \& U_r^\infty(D,D^{\infty}_{d,1}) \& T_r^\infty(D,D^{\infty}_{d,1})    \& L_r^\infty(D,D^{\infty}_{d,1}) \\
    W_r(D,D^{\infty}_{d,1})_\bullet  \& S_r(D,D^{\infty}_{d,1}) \& U_r(D,D^{\infty}_{d,1})\\ 
     \&  \& U_r^{(r-2)}(D,D^{\infty}_{d,1})\& (U_r^{\infty}(D,D^{\infty}_{d,1})^{(r-2)}  \\};
    \path[->]
     (m-1-5) edge node[above]{$\supseteq$}  (m-1-4)
       (m-1-5) edge node[above]{$\cong$}  (m-1-6)
    (m-2-2) edge node[above]{$\pi$}  (m-2-3)
    (m-3-3) edge node[left]{$\subseteq$}  (m-2-3)
     (m-3-3) edge node[above]{$\cong$}  (m-3-4)
       (m-3-4) edge node[left]{$\subseteq$}  (m-1-4)
    (m-2-1) edge node[above]{$\simeq$}  (m-2-2);

   \end{tikzpicture}
  \]
    \caption{A summary of the relationships between the complexes defined so far.}
    \label{fig:sumrelations}
\end{figure}

Let us summarize the relationships we have so far between our various complexes which is illustrated in Figure~\ref{fig:sumrelations}. The leftmost homeomorphism between $W_r(D,D^{\infty}_{d,1})_\bullet$ and $S_r(D,D^{\infty}_{d,1})$ comes from Remark~\ref{rem:firstchange}. By Lemma~\ref{lem-cpt-jn}, there is a complete join map $\pi$ from the complex $S_r(D,D^{\infty}_{d,1})$ to the complex of isotopy classes of asymptotically rigidly embedded loops, $U_r(D,D^{\infty}_{d,1})$, which implies that both complexes have exactly the same connectivity properties. Thus we can choose to work with $U_r(D,D^{\infty}_{d,1})$. Next, Remark~\ref{rem-Un}(1) demonstrates that the $(r-2)$ skeleton of $U_r(D,D^{\infty}_{d,1})$ is the same as the $(r-2)$ skeleton of the complex of asymptotically rigidly embedded loops in $\oplus_r D^{\infty}_{d,1}$, denoted $U_r^{\infty}(D,D^{\infty}_{d,1})$. Since our goal is to show $W_r(D,D^{\infty}_{d,1})_\bullet$ is weakly Cohen-Macauley of dimension $r-2$ (see Corollary~\ref{cor-Sr-conn}), this implies we can again shift our focus to $U_r^{\infty}(D,D^{\infty}_{d,1})$. Next by Definition~\ref{defn-subcpxT}, we have that $T_r^\infty(D,D^{\infty}_{d,1})$ is a subcomplex $U_r^{\infty}(D,D^{\infty}_{d,1})$.  In the next pages, we will show that this complex is isomorphic to the lollipop complex $L_r^\infty(D,D^{\infty}_{d,1})$ as a simplicial complex, Lemma~\ref{lem-id-T-L}, and that the lollipop complex (and hence $T_r^\infty(D,D^{\infty}_{d,1})$) is contractible with a bad simplices argument, Proposition~\ref{prop-Tn-inf}. We then use the contractibility of $T_r^\infty(D,D^{\infty}_{d,1})$ and a bad simplices argument to prove the complex $U_r^{\infty}(D,D^{\infty}_{d,1})$ is contractible and weakly Cohen-Macaulay of dimension $r-2$, Proposition~\ref{prop-Un-inf} and Corollary~\ref{cor-conn-Ur}, implying ultimately that our initial complex is weakly Cohen-Macaulay of dimension $(r-2)$ as needed.

%\end{equation}

%\begin{figure}
 %   \begin{tikzpicture}

 %   \end{tikzpicture}
 %   \caption{My Empty Figure}
  %  \label{compression}
%\end{figure}

%\begin{tikzcd}[scale cd=.85,sep=small]
% & & &  U_r^\infty(D,D^{\infty}_{d,1})  \arrow[hookrightarrow]{r} & T_r^\infty(D,D^{\infty}_{d,1}) &  L_r^\infty(D,D^{\infty}_{d,1}) \\
%W_r(D,D^{\infty}_{d,1})_\bullet   & S_r(D,D^{\infty}_{d,1}) & U_r(D,D^{\infty}_{d,1}) & & & \\
%& & U_r^{(r-2)}(D,D^{\infty}_{d,1})  & U_r^{\infty(r-2)}(D,D^{\infty}_{d,1}) & &
%\end{tikzcd}

In Proposition \ref{prop-Un-inf}, we will deduce the connectivity of $U_r^\infty(D,D^{\infty}_{d,1})$ using the connectivity of the lollipop complex $L_r^\infty(D,D^{\infty}_{d,1})$ by applying a bad simplices argument. Our goal now is to show that $L_r^\infty(D,D^{\infty}_{d,1})$ is highly connected. Let us make some definitions first.

\begin{defn}
Given any lollipop $L: (A,0) \to (D^{\infty}_{d,r},0)$,  we define the \emph{free height} $\mathfrak{h}_L$ to be the minimal  number  $m$ such that $L([1,2])$ is  contained in $D_{d,r,m}$ up to free isotopy.  We also define the height of an admissible loop to be the minimal number $m$ such that it is contained in $D_{d,r,m}$ (see Definition \ref{defn-inf-surf}).
\end{defn}

To analyze the connectivity of $L_r^{\infty}(D,D^{\infty}_{d,1})$, we need the following lemma which is a direct translation of \cite[Lemma 3.8]{SW19}.

\begin{lem}\label{lem-hight-esm}
For any $r,p,N\ge 1$, there exists a number $\mathfrak{h}_{r,p,N}\ge 0$, such that for any $p$-simplex $\sigma$ in $L_r^{\infty}(D,D^{\infty}_{d,1})$, and any $\mathfrak{h}\ge \mathfrak{h}_{r,p,N}$,
there are at least $N$ lollipops of free height $\mathfrak{h}$ in $L_r^{\infty}(D,D^{\infty}_{d,1})$ that are in $\Lk (\sigma)$. 
\end{lem}
\Proof 
Note that for any  vertex $L$ in $L_r^{\infty}(D,D^{\infty}_{d,1})$,  $L\mid_{[1,2]}$ is an admissible loop in  $\oplus_r D^{\infty}_{d,1}$.
Recall the function $q$ defined in Definition \ref{defn-forest} which maps an admissible loop to an edge midpoint in the tree $\CT_{d,r}$. Since each edge has a unique descendent vertex, we can instead map the loop to this vertex which lies in the forest $\CF_{d,r}$. Using this connection, we can now choose $\mathfrak{h}_{r,p,N}$ to be the same as in \cite[Lemma 3.8]{SW19}. Then we have at least $N$ admissible loops of height $\mathfrak{h}\geq \mathfrak{h}_{r,p,N}$ which lie in the complement of the surface corresponding to $\sigma$ in $\oplus_r D^{\infty}_{d,1}$. Connecting each of these admissible loops to the base point in the complement surface, we get a set of lollipops in $\Lk(\sigma)$.
\qed
 
We now show that the complex  $L_r^{\infty}(D,D^{\infty}_{d,1})$ is in fact contractible. The idea of proof is similar to that of  \cite[Proposition 3.1]{SW19} but with significantly more technical difficulty. Intuitively, to define their complex, one only needs information from the loop parts of the lollipops which  are much easier to ``make" them disjoint in general, but for us, we also have to deal with the arc parts which could potentially cause more problems.

\begin{prop}\label{prop-Tn-inf}
The complex $L_r^{\infty}(D,D^{\infty}_{d,1})$ is contractible for any $r\geq 1$.
\end{prop}

\Proof The complex $L_r^{\infty}(D,D^{\infty}_{d,1})$ is obviously non-empty. We will show by induction that for all $k\geq 0$, any map $S^k \rightarrow  L_r^{\infty}(D,D^{\infty}_{d,1})$ is null-homotopic. Assume $L_r^{\infty}(D,D^{\infty}_{d,1})$ is $(k-1)$-connected.

Let $f : S^k \rightarrow  L_r^{\infty}(D,D^{\infty}_{d,1})$ be  a map. As usual, we can assume that the sphere $S^k$ comes with a triangulation such that the map $f$ is simplicial. We first use Lemma \ref{lem:injectifying} to homotope $f$ to a map that is simplexwise injective. For that we need that for every $p$-simplex $\sigma$ in $L_r^{\infty}(D,D^{\infty}_{d,1})$, its link $\Lk{(\sigma)}$ is $(k-p-2)$-connected. But by Lemma \ref{lem-link-loli}, $\Lk({\sigma})$ can be identified with $L_{r_\sigma}^{\infty}(D,D^{\infty}_{d,1})$ for some $r_\sigma\geq 1$, so we have it is $(k-1)$-connected and the conditions of Lemma \ref{lem:injectifying}  are satisfied.

Now since $S^k$ is a finite simplicial complex,  the  free height of the vertices of $S^k$ has a maximum value. We first want to homotope $f$ to a new map such that all the vertices have free height at least $\mathfrak{h}=\mathfrak{h}_{r,k,N}$ where $N =v_0+v_1+\dots+v_k+2$, where $v_i$ is the number of $i$-simplices of $S^k$ and $\mathfrak{h}_{r,k,N}$ is determined by Lemma \ref{lem-hight-esm}. For that we use a bad simplices argument.

We call a simplex of the sphere $S^k$ bad if all of its vertices are mapped to vertices in $ L_r^{\infty}(D,D^{\infty}_{d,1})$ that have free height less than $\mathfrak{h}$.
We will modify $f$ by removing the bad simplices inductively starting by those of the highest dimension. 
Let $\sigma$ be a bad simplex of maximal dimension $p$ among all bad simplices. We will modify $f$ and the triangulation of $S^k$ in the star of $\sigma$ in a way that does not add any new bad simplices. In the process, we will increase the number of vertices by at most $1$ in each step, and not at all if $\sigma$ is a vertex. 
This implies that, after doing this for all bad simplices, we will have increased the number of vertices of the triangulation of $S^k$ by at most $v_1+\dots+v_k$. As $S^k$ originally had $v_0$ vertices, at the end of the process its new triangulation will  have at most $v=v_0+v_1+\dots+v_k$ vertices. There are two cases.

\textbf{Case 1: $p=k$.} If the bad simplex $\sigma$ is of the dimension $k$ of the sphere $S^k$, then its image $f(\sigma)$ has a complement loop which bounds a surface $C$ asymptotically rigidly homeomorphic to $D^{\infty}_{d,r_\sigma}$ for some $r_\sigma\geq 1$ by Lemma \ref{lem-link-loli}. Now we can choose a lollipop $y$ in $C$ with free height at least $\mathfrak{h}+1$. In particular $f(\sigma) \cup y$ form a $(k+1)$-simplex. We can then add a vertex $a$ in the center of $\sigma$, replacing $\sigma$ by $\partial \sigma \ast a$ and replacing $f$ by the map $(f|_{\partial\sigma})\ast (a \mapsto y)$ on $\partial\sigma \ast a$. This map is homotopic to $f$ through the simplex $f(\sigma)\cup \{y\}$. We have added a single vertex to the triangulation. Because $L$ has free height $\mathfrak{h}+1$, we have not added any new bad simplices, and we have removed one bad simplex, namely $\sigma$. Moreover, $f$ remains simplexwise injective. 

\textbf{Case 2: $p<k$.} If the bad simplex $\sigma$ is a $p$-simplex for some $p<k$, by maximality of its dimension, the link of $\sigma$ is mapped to vertices of free height at least $\mathfrak{h}$ in the complement of the subsurface $f(\sigma)$. 
The simplex $\sigma$ has $p+1$ vertices whose images are pairwise disjoint outside the based point up to based isotopy.
By Lemma \ref{lem-hight-esm} and  our choice of $\mathfrak{h}$, there are at least $N= v+2$ lollipops  $y_1,\dots,y_N$ of free height $\mathfrak{h}$ such that each $f(\sigma)\cup\{y_i\}$ form a $(p+1)$-simplex. As there are fewer vertices in the link than in the whole sphere $S^k$, and $S^k$ has at most $v$ vertices, by the pigeonhole principle, the loop part of the vertices in $f(\Lk(\sigma))$ are contained in at most $v$ punctured disks bounded by the corresponding admissible loops with free height $\mathfrak{h}$. As $N=v+2$, 
there are at least two of the above vertices $y_i$ and $y_j$ of free height $\mathfrak{h}$ such that any loop parts of  vertices in $f(\Lk(\sigma))$ are disjoint from the loop parts of $y_i$ and $y_j$. We can further assume that the arc parts of $y_i$ and $y_j$ never intersect with any loop part of the vertices in $f(\Lk(\sigma))$. And up to replacing the loop part of $y_i$ and $y_j$ with an admissible loop lying inside the disk bounded by the loop parts of  $y_i$ and $y_j$ (note that this may increase the free height of $y_i$ and $y_j$), we can further assume that the arc parts of  vertices in $f(\Lk(\sigma))$ are disjoint from the loop part of $y_i$ and $y_j$.  But unlike the situation in the proof of \cite[Proposition 3.1]{SW19}, a new  problem we are facing here is that the arc parts of $y_i$ or $y_j$ might intersect the arc parts of the vertices in $f(\Lk(\sigma))$ even up to  isotopy. In particular,  given a simplex $\tau$ lying in the link of $\sigma$,  $f(\sigma)\cup f(\tau)\cup y_i$ does not necessarily form a simplex now.

 For that we want to apply  the  mutual link trick (cf. Lemma \ref{lemma-replace-trick}) to remove the intersections of $f(\Lk(\sigma))$ with $y_i$ via a sequence of homotopies. In the process, we will only modify $f$ on $\Lk(\sigma)$ and the new map still maps $\Lk(\sigma)$ to $\Lk_{L_r^{\infty}(D,D^{\infty}_{d,1})} (f(\sigma))$. Recall that $f$ is simplexwise injective.   Up to isotopy, we can further choose representatives for vertices in $f(\Lk(\sigma))$ such that the intersection points of vertices in $f(\Lk(\sigma))$ and $y_i$ are isolated. Moreover, we assume the number of intersection points is minimal for each vertex in $f(\Lk(\sigma))$. Now we choose an intersection point $x_0$ in the arc $y_i([0,1])$ that is closest to $y_i(1)$, denote the corresponding lollipop by $\beta$ which is the image of some vertex $b \in \Lk (\sigma)$. We can choose $\beta'$ to be a  variation of $\beta$: $\beta'$ coincides with $\beta$ for the most part, except around the intersection point with $y_i$, we replace it by an arc going around the loop part of $y_i$. See Figure \ref{jumping} for a picture of this. Now we apply Lemma \ref{lemma-replace-trick}, for which we need to check the following two conditions: 
 
 \begin{figure}
\centering
\includegraphics[width=.6\textwidth]{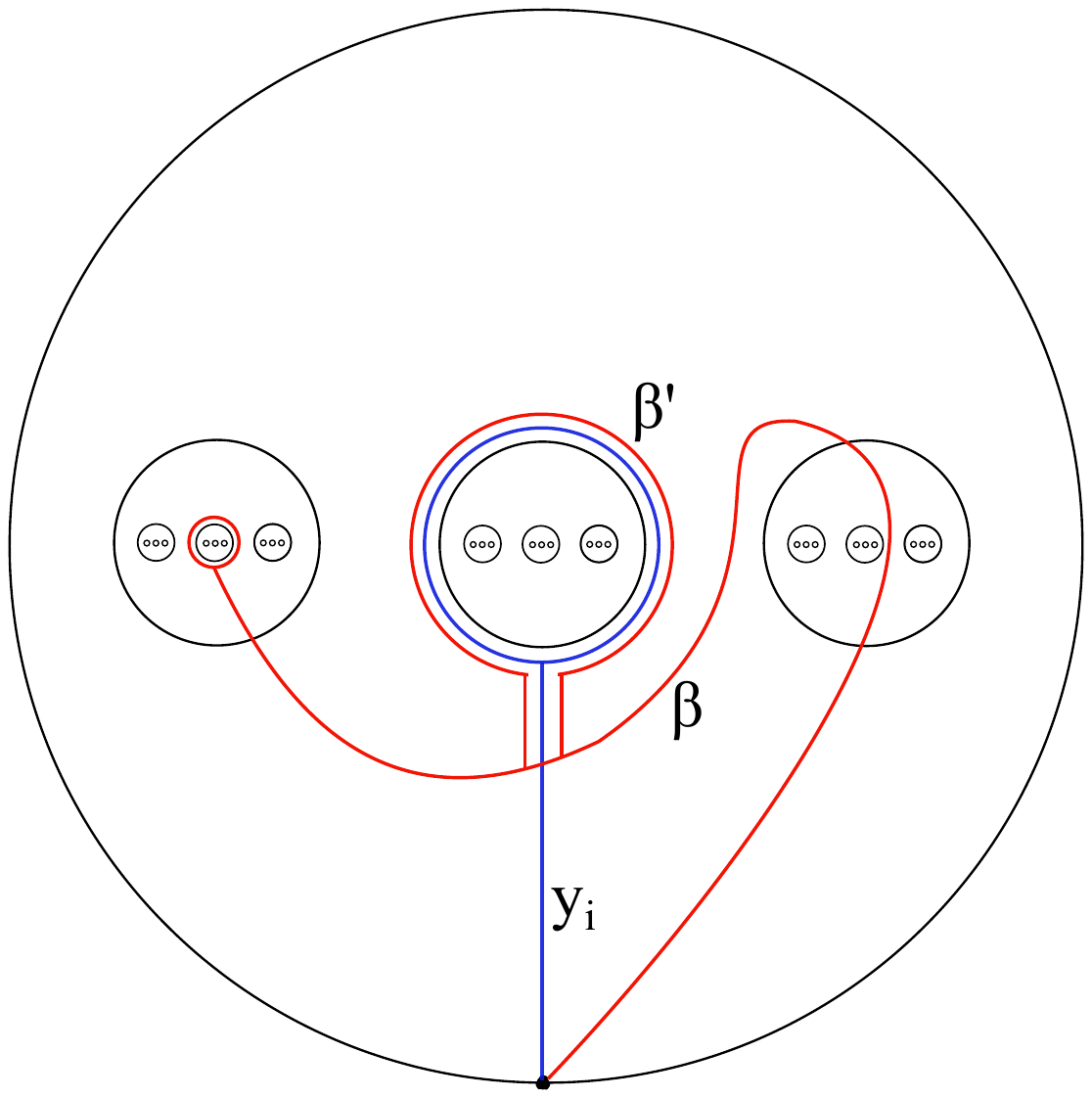}
\caption{Replacing $\beta$ by $\beta'$ to reduce the number of intersection points with $y_i$.}
\label{jumping}
\end{figure}

 \begin{enumerate}
     \item $f(\Lk_{S^k}(b)) \leq \Lk_{L_r^{\infty}(D,D^{\infty}_{d,1})} (\beta')$. This follows from our definition of $\beta'$. If a vertex $v$ in $f(\Lk_{S^k}(b))$ is disjoint from $\beta$, using the fact that the intersection point $x_0$ is the closest one  to $y_i(1)$ and $f(v)$ is disjoint from the loop part of $y_i$, we have $\beta'$ is also disjoint from $v$.
     \item  $\Lk_{} (\beta) \cap \Lk_{}(\beta')$ is $(k-1)$-connected. The lollipops $\beta$ and $\beta'$ together will bound a disk which contains the loop part of  $\beta$ and  $y_i$. In any event, the complement of these is a surface asymptotically rigidly homeomorphic to some surface $D^{\infty}_{d,r'}$ for some $r'\geq 1$. By our induction, it is $(k-1)$-connected. 
 \end{enumerate}
Now  Lemma \ref{lemma-replace-trick} says we can homotope $f$ to a new map such that $f(b) =\beta'$ and $f(\Lk(\sigma))$ has fewer intersection points with $y_i$. Step by step, at the end we have a simplexwise injective map $f$ such that for any vertex in $f(\Lk(\sigma))$, it only intersects with $y_i$ at the base point. In particular for any $\tau \in \Lk(\sigma)$, we have $f(\sigma)\cup f(\tau)\cup\{y_i\}$ forms a simplex in $L_r^{\infty}(D,D^{\infty}_{d,1})$. 

%Now each $f(\tau)\cup f(\sigma)\cup\{y_i\}$ forms a simplex of $L_r^{\infty}(D,D^{\infty}_{d,1})$. 
We can then replace $f$ inside the star
\[
\St(\sigma)=\Lk(\sigma)*\sigma \simeq S^{k-p-1}*D^p
\]
by the map $(f|_{\Lk(\sigma)})*(a\mapsto y_i)*(f|_{\partial \sigma})$ on 
\[
\Lk(\sigma)*a*\partial\sigma \simeq S^{k-p-1}*D^0*S^{p-1}.
\]
which agrees with $f$ on the boundary $\Lk(\sigma)*\partial\sigma$ of the star, and is homotopic to $f$ 
through the map $(f|_{\Lk(\sigma)})*(a\mapsto y_i)*(f|_{\sigma})$ defined on 
\[
\Lk(\sigma)*a*\sigma\simeq S^{k-p-1}*D^0*D^p.
\] 
Now $\Lk(\sigma)*a*\partial(\sigma)$  has exactly one extra vertex $a$ compared to the star of $\sigma$, unless $\sigma$ is just a vertex, in which case its boundary is empty and it has the same number of vertices.
As $y_i$ has height at least $\mathfrak{h}$, we have not added any new bad simplices. Hence we have reduced the number of bad simplices by one by removing $\sigma$. 

By induction, we can now assume that there are no bad simplices for $f$ with respect to a triangulation with at most $v$ vertices. With this assumption, we want to 
cone off $f$ just as we coned off the links in the above argument.  We have more than $N=v+2$ vertices of free height $\mathfrak{h}$ in $L_r^{\infty}(D,D^{\infty}_{d,1})$, and at most $v$ vertices in the sphere. The loop parts of these vertices are admissible loops of  height at least $\mathfrak{h}$. By the pigeonhole principle, we know that there are at least two lollipops $z_i$ and $z_j$ of free height $\mathfrak{h}$ such that the punctured disks bounded by their loop parts are disjoint from the punctured disk bounded by any loop part  of the lollipops in the vertices of $f(S^k)$. Just as before, we can further assume that the arc parts of $z_i$ and $z_j$ never intersect with any loop part of the vertices in $f(S^k)$, and the arc parts of  vertices in $f(S^k)$ are disjoint from the loop part of $z_i$ and $z_j$. But the same problem appears again, as we want vertices of $f(S^k)$ to be disjoint from the whole lollipop $z_i$. For that we apply  Lemma \ref{lemma-replace-trick} again and the same proof as before implies that we can homotope $f$ such that its image is disjoint from $z_i$. In particular $f(S^k)$ lies in the link of $z_i$. Hence we can homotope $f$ to a constant map since $\St(z_i)$ is contractible.
\qed

\begin{prop}\label{prop-Un-inf}
The complex $U_r^{\infty}(D,D^{\infty}_{d,1})$ is contractible.
\end{prop}

\Proof As $T_r^{\infty}(D,D^{\infty}_{d,1})$ is a subcomplex of $U_r^{\infty}(D,D^{\infty}_{d,1})$, we can use a bad simplices argument. 

We call a vertex of $U_r^{\infty}(D,D^{\infty}_{d,1})$ bad if it does not lie in $T_r^{\infty}(D,D^{\infty}_{d,1})$ and a simplex bad if all of its vertices are bad.  Given a bad $p$-simplex $\sigma$, we need to determine the connectivity of the good link $GL_{\sigma}$ (see Subsection \ref{sub-bad-sim} for the definition of $GL_\sigma$).   As in the proof of Lemma \ref{lem-link-loli}, we have a complement surface $C_\sigma$ of $\sigma$ in $D^{\infty}_{d,1}$. Note that $C_\sigma$ inherits a $d$-rigid structure and  it is asymptotically rigidly homeomorphic to $\oplus_{r_\sigma}D^{\infty}_{d,1} $ for some $r_\sigma>0$.  In particular, we can now identify $GL_\sigma$ with $T_{r_\sigma}^{\infty}(D,D^{\infty}_{d,1})$ which is contractible. Thus by Proposition \ref{prop-bad-sim}, we have the pair $(U_r^{\infty}(D,D^{\infty}_{d,1}),T_r^{\infty}(D,D^{\infty}_{d,1}))$ is $i$-connected for any $i\geq 0$. By Proposition \ref{prop-Tn-inf}, $T_r^{\infty}(D,D^{\infty}_{d,1}) \cong L_r^{\infty}(D,D^{\infty}_{d,1})$ is contractible, so we also have $U_r^{\infty}(D,D^{\infty}_{d,1}))$ is contractible.
\qed

\begin{cor}\label{cor-conn-Ur}
The complex $U_r(D,D^{\infty}_{d,1})$ is weakly Cohen-Macaulay of dimension $r-2$.
\end{cor}

\Proof 
Note first that a simplicial complex is $(r-3)$-connected if and only if its $(r-2)$-skeleton is. Since $U_r(D,D^{\infty}_{d,1})$ has the same $(r-2)$-skeleton as $U_r^{\infty}(D,D^{\infty}_{d,1})$ and $U_r^{\infty}(D,D^{\infty}_{d,1})$ is contractible, in particular $(r-3)$-connected, we indeed have $U_r(D,D^{\infty}_{d,1})$ is $(r-3)$-connected.  

Now let $\sigma$ be a $p$-simplex of $U_r(D,D^{\infty}_{d,1})$, with vertices $\phi_0,\phi_1,\cdots,\phi_p$. We need to check that the link  $\Lk_ {U_r(D,D^{\infty}_{d,1})}(\sigma)$ is $(r-p-4)$-connected. We can assume $p\leq r-3$ as any space is $(-2)$-connected. Moreover, it suffices to show the $(r-p-3)$-skeleton of $\Lk_{ U_r(D,D^{\infty}_{d,1})}(\sigma)$ is $(r-p-4)$-connected. Since $\phi_0,\phi_1,\cdots,\phi_p$ forms a $p$-simplex, similar to the proof of Lemma \ref{lem-link-loli}, we have the complement surface of $\sigma$ is asymptotically rigidly homeomorphic to some $d$-rigid surface $D^{\infty}_{d,k_\sigma}$ for some $k_\sigma>0$. Then  we can identify the  $(r-p-3)$-skeleton of $\Lk_{U_r(D,D^{\infty}_{d,1})}(\sigma)$ with the  $(r-p-3)$-skeleton of $U_{k_{\sigma}}^{\infty}(D,D^{\infty}_{d,1})$. Since  $U_{k_\sigma}^{\infty}(D,D^{\infty}_{d,1})$ is even contractible, we have the connectivity bound we need.
\qed

Now by Lemma \ref{lem-cpt-jn} and Proposition \ref{prop-join-conn}, we have the following.

\begin{cor}\label{cor-Sr-conn}
The complexes $S_r(D,D^{\infty}_{d,1})$ and $W_r(D,D^{\infty}_{d,1})_\bullet$ are weakly Cohen-Macaulay of dimension $r-2$.
\end{cor}

\subsection{Homological stability} We are finally ready to prove the homological stability result.

\begin{thm}\label{thm-hmg-stab-RB+}
Suppose $d\geq 2$. Then the inclusion maps induce isomorphisms
$$ \iota_{R^+,d,r}: H_i(RV^+_{d,r},M) \to H_i(RV^+_{d,r+1},M)$$
in homology in all dimensions $i \geq 0$, for all $r \geq  1$ and for all $H_1(RV^+_{d,\infty})$-modules $M$.
\end{thm}

\Proof From Corollary \ref{cor-Sr-conn},
we have that $W_r(D, D_{d,1}^\infty)$ is $(r - 2)$-connected, hence in particular, the category $\mathcal{G}_d$ satisfies property $LH3$ at the pair of objects $(D, D_{d,1}^\infty)$ with slope k = 3. By Theorem \ref{thm:hom stab},  we have for any abelian $RV^+_{\infty}$-module $M$ the map
\[H_i(RV^+_{d,r}; M )  \longrightarrow H_i(RV^+_{d,r+1}; M)\]
induced by the natural inclusion map is isomorphism if   $r \geq 3i+3$.

But we can improve the stability range as in the proof of \cite[Theorem 3.6]{SW19} by noticing that we have the same canonical isomorphism between $RV^+_{d,r}$ and $RV^+_{d,r+d-1}$. In fact, denote this isomorphism by $I_{d,r}$, we have the following commutative diagram:

\[\begin{tikzcd}
RV^+_{d,1+ r-1 }\arrow[d,"I_{d,r}"] \arrow[r, "\iota_{R^+,d,r}"] & ~ RV^+_{d,1+(r-1) +1}\arrow[d,"I_{d,{1+(r-1) +1}}"] \\
RV^+_{d,{d+{r-1}}}  \arrow[r, "\iota_{R^+,d,r +1}"] & ~RV^+_{d,{d+{r-1}+1}}
\end{tikzcd}
\]

Given that the vertical maps are isomorphism, and the bottom horizontal maps induce isomorphism on the $i$-th homology when   $d+r-1 \geq 3i+3$, we also the top map must also induce isomorphism on the homology as long as   $r \geq 3i+3$. This has improved the stable range by $d-1$. Step by step, we must have the map $\iota_{R^+,d,r}$ induce isomorphism on homology in  all dimensions $i\geq 0$ and for all $r\geq 1$. 

This finishes the proof of Theorem \ref{thm-hmg-stab-RB+}.
\qed

\begin{thm}\label{thm-hmg-stab-RB}
Suppose $d\geq 2$. Then the inclusion maps induce isomorphisms
$$ \iota_{R,d,r}: H_i(RV_{d,r},M) \to H_i(RV_{d,r+1},M)$$
in homology in all dimensions $i \geq 0$, for all $r \geq  1$ and for all $H_1(RV_{d,\infty})$-modules $M$.
\end{thm}
\sProof The proof will be exactly the same as that of Theorem \ref{thm-hmg-stab-RB+}.  Note first that by Theorem \ref{thm-iden-asym-ribb}, it is the same as proving the half-twist asymptotic mapping class groups $ \CH V_{d,r}$  have homological stability.  We define the braided monoidal category $\CG'_d$ to be the category with objects $\oplus_r D^{\infty}_{d,1}$, $r\ge 0$, $\oplus$ as the operation, and $D$ as the $0$ object. When $r=s$, we define the morphisms $\Hom (\oplus_r D^{\infty}_{d,1}, \oplus_s D^{\infty}_{d,1}) = \CH V_{d,r}$  which can also be understood as the group of isotopy classes of  asymptotically quasi-rigid homeomorphisms of $\oplus_r D^{\infty}_{d,1}$; when $r\neq s$, let $\Hom (\oplus_r D^{\infty}_{d,1}, \oplus_s D^{\infty}_{d,1}) =\emptyset $. We then have a homogeneous category $U\CG'_d$ and to prove the homological stability for the sequence of  groups $\CH V_{d,1}\leq \CH V_{d,2}\leq \cdots$, we only need to prove the associated space $W_r(D,D^{\infty}_{d,1})_\bullet$, in fact the associated simplicial complex $S_r(D,D^{\infty}_{d,1})$, is highly connected. At this point, the complex is slightly different from the oriented case, but still the new complex $S_r(D,D^{\infty}_{d,1})$ is a complete join over the old complex $U_r(D,D^{\infty}_{d,1})$. Hence, the connectivity of $S_r(D,D^{\infty}_{d,1})$ again follows from Corollary \ref{cor-conn-Ur} and Proposition \ref{prop-join-conn}.
\qed

\bibliographystyle{alpha}
\bibliography{references.bib}

\end{document}